\documentclass{amsart}
\usepackage{amssymb}

\usepackage[all,cmtip]{xy}

\usepackage[usenames,dvipsnames]{color}
\usepackage{hyperref}

\usepackage{graphicx}

\usepackage[ansinew]{inputenc}

   {\end{trivlist}\platz}

\newtheorem{theorem}{Theorem}[section]
\newtheorem{proposition}[theorem]{Proposition}
\newtheorem{lemma}[theorem]{Lemma}

\newtheorem{definition}[theorem]{Definition}

\newtheorem{remark}[theorem]{Remark}

\newcommand {\eps}{\varepsilon}
\newcommand\calC{{\mathcal C}}

\newcommand\calK{{\mathcal K}}

\newcommand\htilde{\tilde{h}}

\newcommand{\Etilde}{\widetilde{E}}
\newcommand{\Gtilde}{\widetilde{G}}
\newcommand{\Wtilde}{\widetilde{\bW}}
\newcommand{\Xtilde}{{\widetilde{X}}}
\newcommand{\Ztilde}{{\widetilde{Z}}}
\newcommand\Id{{\rm Id}}

\newcommand {\R}{\mathbb{R}}

\newcommand {\N}{\mathbb{N}}

\newcommand {\Q}{\mathbb{Q}}

\newcommand{\bg}{{\bf g}}

\newcommand{\bo}{{\bf 0}}
\newcommand{\bp}{{\bf p}}

\newcommand{\bu}{{\bf u}}
\newcommand{\bq}{{\bf q}}
\newcommand{\bS}{{\bf S}}
\newcommand{\bV}{{\bf V}}

\newcommand{\bW}{{\bf W}}

\newcommand{\dpl}{\partial}

\newcommand{\partildX}{{\partial \Xtilde}}
\newcommand{\tildX}{\Xtilde}  %% unclear, maybe widetilde?

\newcommand{\Sr}{{\bf S}}

\newcommand{\vp}{\varphi}

\newcommand{\overbp}{{\overline{\bf p}}}
\newcommand{\overbq}{{\overline{\bf q}}}

\newcommand{\rdot}{\dot{r}}

\newcommand{\xidot}{\dot{\xi}}

\newcommand{\thetadot}{\dot{\theta}}

\newcommand{\geucl}{g_{\text{eucl}}}

\newcommand{\threeT}{{}^{\text{3}}T}
\newcommand{\Sbound}{S}  %% changed from S^\partial in Aug 2011
\newcommand{\vpdot}{\dot{\vp}}

\newcommand{\ftilde}{\tilde{f}}

\newcommand{\Ebound}{E^\partial}
\newcommand{\Ebounddot}{\dot{E}^\partial}

\newcommand{\gdM}{\bg_{\partial M}}
\newcommand{\gdXtilde}{\bg_{\partial \Xtilde}}

\newcommand{\bpmin}{{\bp_{\text{min}}}}
\newcommand{\bpmax}{{\bp_{\text{max}}}}
\newcommand{\gammamin}{{\gamma_{\text{min}}}}
\newcommand{\overbpmin}{\overbp_{\text{min}}}
\newcommand{\Gammabpm}{\Gamma_{\bp_-}}
\newcommand{\vpmin}{\vp_{\text{min}}}
\newcommand{\interior}[1]{\mathring{#1}}
\newcommand{\dM}{{\partial M}}
\newcommand{\intM}{\interior{M}}
\newcommand{\expdM}{\exp_{\dM}}

\newcommand{\Cpartial}{C_\partial}  % C restricted to boundary
\newcommand{\Utilde}{\tilde U}
\begin{document}
\title{The exponential map at a cuspidal singularity}
\author{Vincent Grandjean}
\address{Departamento de Matem\'atica, UFC,
Av. Humberto Monte s/n, Campus do Pici Bloco 914,
CEP 60.455-760, Fortaleza-CE, Brasil}
\email{vgrandje@mat.ufc.br}
\author{Daniel Grieser}
\address{Institut f\"ur Mathematik, Carl von Ossietzky Universit\"at Oldenburg, 26111 Oldenburg, Germany}
\email{daniel.grieser@uni-oldenburg.de}
\thanks{The first author is very grateful to the Fields Institute for working conditions and support 
while completing this work. He is also very pleased to thank Carl von Ossietzky Universit\"at Oldenburg for support while 
visiting twice the second author.
Both authors were supported in part by the Deutsche Forschungsgemeinschaft in the priority program 
'Global Differential Geometry'.
Both authors are very grateful to Vincent Naudot for his help and advice on linearization.}
\subjclass[2010]{
53B21   	%Methods of Riemannian geometry
37C10   	%Vector fields, flows, ordinary differential equations
%53C50   	Lorentz manifolds, manifolds with indefinite metrics (53C is global diff. geometry)
53C22   	%Geodesics [See also 58E10]
%(34: ODEs 34C: Qualitative theory)
%34F15   	Resonance phenomena
%34C05   	Location of integral curves, singular points, limit cycles
%34C45   	Invariant manifolds
%(32: Several complex var; 32S: singularities
%32S70   	Other operations on singularities
%(58: Global Analysis; 58K: Theory of singularities)
%58K45   	Singularities of vector fields, topological aspects
%53B20%Local Riemannian geometry
}

\begin{abstract}
We study spaces with a cuspidal (or horn-like) singularity embedded in a smooth Riemannian manifold and analyze 
the geodesics in these spaces which start at the singularity. 
This provides a basis for understanding the intrinsic geometry of such spaces near the singularity.
We show that these geodesics combine to naturally define an 
exponential map based at the 
singularity, but that the behavior of this map can deviate strongly 
from the behavior of the exponential map based at a smooth point or at a conical singularity: While it 
is always surjective near the singularity, it may be discontinuous and non-injective on any neighborhood of the 
singularity. The precise behavior of the exponential map is determined by a function on the link of the singularity which is an 
invariant of the induced metric.
Our methods are based on the Hamiltonian system of geodesic differential equations and on techniques of singular analysis.
 The results are proved
in the more general natural setting of manifolds with boundary carrying a so-called cuspidal metric.
\end{abstract}
\date{\today}
\maketitle
\setcounter{tocdepth}{1} % subsections don't appear in table of contents
\tableofcontents

\section{Introduction}
Geodesics are among the fundamental objects of differential geometry.
On a smooth Riemannian manifold the geodesics starting at a point $p$ are classified by and smoothly 
depend on their initial velocity vector, and combine to define the exponential map based at 
$p$, which in turn yields normal coordinates and important geometric information about the manifold. Also, 
geodesics arise in the study of the propagation of waves on the manifold as paths along which 
singularities of solutions of the linear wave equation travel.

Much less is known about the behavior of geodesics on 
{\em singular} spaces. Here by a geodesic we always mean a locally shortest curve.  Our general aim is to analyze the full local asymptotic behavior of the family 
of geodesics reaching, or starting at, a singularity. Previously, Bernig and Lytchak \cite{BerLyt:TSGHLSS} 
obtained first order information for geodesics on general real algebraic sets  $X\subset \R^n$, by 
showing that any geodesic reaching a singular 
point $p$ of $X$ in finite time must have  a limit direction at $p$.
In the special situation of isolated conical singularities $p$,
Melrose and Wunsch \cite{MelWun} obtain full asymptotic information by showing that the geodesics 
starting at $p$ define a smooth foliation of a neighborhood of $p$ 
and thus may be combined to define a smooth exponential map based at $p$, analogous to Theorem \ref{thm:S=const} below.

In this paper we analyze the family of geodesics starting at an isolated cuspidal singularity, defined below. Our method is based on the geodesic differential equations and yields full asymptotic information about the geodesics.
We will see that there is a much richer range of possible local behavior
than in the case of conical singularities; for example, there is a natural notion of exponential map based 
at the singularity, but this map may be non-injective or discontinuous on any neighborhood of the singularity. 

A natural setting for our investigation is the notion of {\em cuspidal manifold}. This is a 
smooth (that is, $C^\infty$) manifold, $M$, with compact boundary, equipped with a semi-Riemannian
 metric, $\bg$, which is Riemannian in the interior $\interior{M}$ and in a neighborhood 
$U$ of the boundary $\partial M$ can be written
\begin{equation}
 \label{eqn:cusp metric inv k}
 \bg_{|U} = \left(1-k(k-1)r^{2k-2} S + O(r^{2k-1})\right)\,dr^2 + r^{2k} h \,.
\end{equation}
Here $r$ is a boundary defining function for $M$ (that is, $r\in C^\infty(M)$ is positive 
in $\interior{M}$, vanishes on $\partial M$ and satisfies $dr_{\bp}\neq 0$ for all 
$\bp\in\partial M$), $S$ is a smooth function on $\partial M$ and $h$ is a smooth symmetric 
two-tensor on $U$ whose restriction to $\partial M$, denoted $\gdM$, is positive definite. 
Finally, 
$$k\geq2\text { is an integer, called the {\em order} of the cuspidal manifold.}$$
At first reading it may be useful to think of $k=2$, so that  $\bg_{|U} = (1-2r^2 S + O(r^3))\,dr^2 + r^4 h$.
In fact, it follows from our results that cusps of different orders behave qualitatively the same with respect to the questions we investigate.

We 
call $\bg$ a {\em cuspidal metric of order $k$}. 
A given 
cuspidal metric defines invariantly the quantities 
$$ S\in C^\infty(\partial M)\quad\text{and}\quad \gdM, \text{ a Riemannian metric on }\partial M$$
(more precisely, $S$ is only determined up to addition of a locally constant function)
and fixes $r$ to order $2k-1$. For simplicity we will fix a boundary defining function $r$ throughout.

\begin{figure}
 \includegraphics[width=6cm]{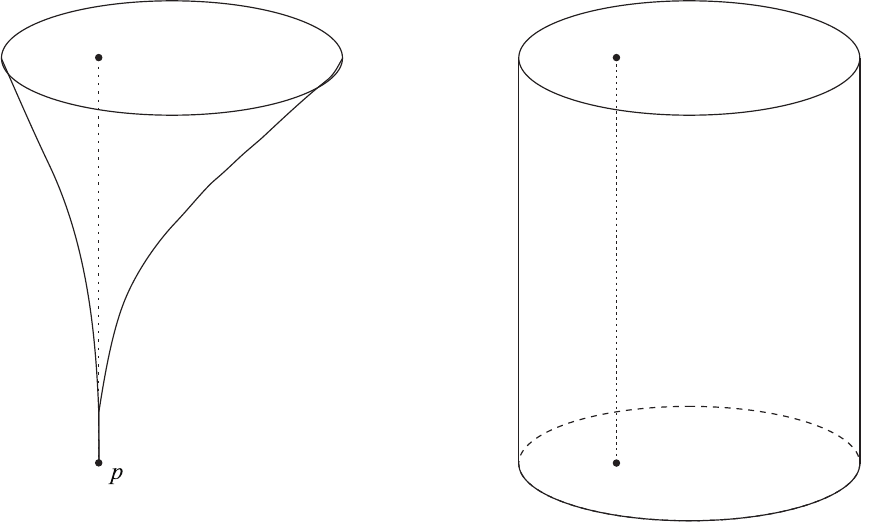}
 \caption{A cuspidal surface singularity $X$ and its resolution $\Xtilde$; the dotted lines are explained in Remark \ref{rem:S, gdm for Xtilde}}
\label{fig:cusp}
\end{figure}

Cuspidal manifolds of order $k$ arise as resolutions of spaces $X$ with isolated cuspidal singularity $p\in X$ or order $k$. These are subsets $X\subset\R^n$ so that $X\setminus p$ is a submanifold 
of $\R^n$ and so that $X$ is given by the 
following local model in a neighborhood $U'$ of $p$: Let $p=0$. Then, with $\R_+=[0,\infty)$,
\begin{equation}
 \label{eqn:beta}
X\cap U' = \beta (\tildX),\quad \beta : \R^{n-1}\times \R_+ \to \R^{n},\quad (u,z) \mapsto (z^ku,z)
\end{equation}
where $\tildX\subset\R^{n-1}\times \R_+$ is a smooth manifold with compact boundary, which satisfies
\begin{equation}
\label{eqn:beta'}
 \tag{\ref{eqn:beta}'}
\partial \tildX \subset E:=\R^{n-1}\times \{0\},\quad \Xtilde\text{ intersects $E$ transversally.} 
\end{equation}
We call 
$\Xtilde$ a {\em resolution} of $X$. See Figure \ref{fig:cusp}.
Now if $g$ is any smooth Riemannian metric on the ambient space $\R^n$ then the restriction of $g$ to
the manifold $X\setminus p$ pulls back under $\beta$ to a Riemannian metric on the interior 
of $\Xtilde$.
We will prove that this metric extends  smoothly to the boundary, and that the extension is a cuspidal metric on $\Xtilde$. 
In this sense the resolution of a space with cuspidal singularity is a cuspidal manifold $(M,\bg)=(\Xtilde,\beta^*g)$.

Let $(M,\bg)$ be a cuspidal manifold. Since $\bg$ is Riemannian over the interior $\intM$, locally shortest curves in $\intM$ 
are precisely the unit speed 
solutions of the geodesic differential equations. We call these solutions {\em geodesics on $M$}.  We are interested 
in those geodesics {\em starting at the boundary}, that is, geodesics $\gamma:(T_-,T_+)\to\interior{M}$ 
for which $\lim_{\tau\to T_-} r(\gamma(\tau)) = 0$. Here $T_-\in\R\cup\{-\infty\}$. If $\lim_{\tau\to T_-}\gamma(\tau)$ exists
 then we call it the {\em starting point} of $\gamma$.

Our first theorem shows that, under certain assumptions, geodesics starting at the 
boundary have a starting point on $\partial M$, and characterizes the possible starting points. 
Recall that a critical point of $S\in C^\infty(\partial M)$ is a point $\bp\in\partial M$ where 
the differential $dS_{|\bp}$ vanishes.
\begin{theorem} \label{thm: crit pts}
Let $(M,\bg)$ be a cuspidal manifold.
Assume that either $S$ is constant or the critical points of $S$ are isolated. 

Then there are $r_0>0$, $\tau_0>0$ so that the following is true:
Let $\gamma:(T_-,T_+)\to\interior{M}$ be a maximal geodesic and assume that for 
some $T\in(T_-,T_+)$  the initial part $\gamma_{|(T_-,T)}$ is contained 
in $\{r<r_0\}$. Then
\begin{enumerate}
 \item[(a)]
 $T_-$ is finite and the starting point $\lim_{\tau\to T_-}\gamma(\tau)$ exists. 
 \item[(b)] The starting point of $\gamma$ lies on $\partial M$ and is a critical point of $S$.
 \item[(c)] $T_+-T_- \geq \tau_0$, i.e. the geodesic exists at least for time $\tau_0$. 
\end{enumerate}
Also, the following converse to (b) holds: Each critical point of $S$ is the starting point of a geodesic on $M$. 
\end{theorem}

Our next theorem says that for constant $S$ we have a smooth foliation by 
geodesics of a neighborhood of $p$,  analogous to the Theorem of Melrose and Wunsch 
\cite{MelWun} in the conical case.

\begin{figure}
  \includegraphics[width=4cm]{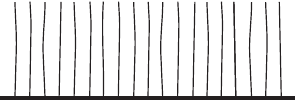}
  \includegraphics[width=4cm]{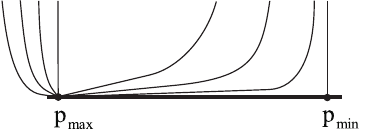}
  \includegraphics[width=4cm]{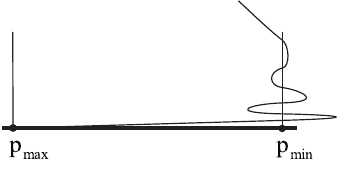}
 \caption{Geodesics in the cases $S=$ const, $S_{\vp\vp}<a_k$ and $S_{\vp\vp}(\bpmin)>a_k$. 
The bold line is $\dM$. The points $\bpmax$ and $\bpmin$ are  a maximum resp.\ minimum of $S$.}
 \label{fig:thms}
\end{figure}

\begin{theorem}
\label{thm:S=const}
Let $(M,\bg)$ be a cuspidal manifold.
Assume $S$ is constant. 
 Then to each $\bp\in\partial M$ there is a unique geodesic $\gamma_\bp$ starting at $\bp$ at time $\tau=0$.
Furthermore, there is $\tau_0>0$ and a neighborhood $U$ of $\partial M$ in $M$ such that the exponential map
$$ \expdM: \partial M \times [0,\tau_0) \to U\subset M,\quad (\bp,\tau) \mapsto \gamma_\bp(\tau)$$
is a diffeomorphism.
\end{theorem}
See the left picture in Figure \ref{fig:thms}.
For a general cuspidal manifold, where $S$ need not be constant, it seems clear that a full neighborhood 
of the boundary is covered by geodesics starting at the boundary, since for any interior point of $M$ 
any shortest curve from this point to the boundary must be a geodesic. We will not attempt to make this 
argument precise here (that is, prove existence of a minimizer) but rather analyze the case where $S$ is a Morse function much more precisely, 
mostly in the case of surfaces.

For our next theorems we assume that $S$ is a Morse function. That is, all critical points are 
non-degenerate (i.e.\ the Hessian of $S$ at these points is a non-degenerate quadratic form). 
This implies that the critical points are isolated.

The behavior of geodesics hitting a non-degenerate critical point $\bp$ of $S$ can be analyzed 
very precisely in a neighborhood of $\bp$. See for example Proposition \ref{prop-localmax} for 
the case of local maxima and minima of $S$. We will now focus on surfaces $M$, for which we can 
analyze how these local pictures fit together in a full neighborhood of $\partial M$. 
First, we show that there is an exponential map.
To simplify the exposition we assume in the following theorems that $\dM$ is connected.
 
 \begin{theorem}
 \label{thm:surfaces}
 Let $(M,\bg)$ be a cuspidal surface with connected boundary. Assume that $S$ is a Morse function. 
There is a parametrization $q\mapsto \gamma_q$, $q\in \bS^1$, of the set of geodesics starting at $\dM$ 
at time $\tau=0$ so that the map
\begin{equation}
 \label{eqn:exp map general}
\expdM : \bS^1\times (0,\tau_0) \to U\setminus\dM,\quad (q,\tau)\mapsto \gamma_q(\tau)
\end{equation}
is defined and surjective for some $\tau_0>0$ and some neighborhood $U$ of $\dM$. 

The full asymptotic behavior of $\expdM$ as $\tau\to 0$ can be described explicitly in terms of $S$.
\end{theorem}

The parametrization $q\mapsto\gamma_q$ is uniquely determined, up to homeomorphisms of $\bS^1$, by 
the requirement that it preserves cyclic ordering in a suitable sense, explained after Definition \ref{def:exp map}. However, $\expdM$ need not be 
a diffeomorphism, it may even be discontinuous. See Section \ref{sec:exponential map} for details.

Whether the exponential map is a homeomorphism onto a neighborhood of the boundary or not is mostly determined by the size of the second derivative of $S$, and the number
\begin{equation}
\label{eqn:def ak} 
 a_k = \frac{(2k-1)^2}{2 k(k-1)} = 2 + \frac 1{2k(k-1)}
\end{equation}
(for example, $a_2 = 9/4$) 
turns out to be the determining threshold, where $k$ is the order of the cusp. 
This is made precise in the next two theorems. Here 
and throughout we denote by $S_\vp$, $S_{\vp\vp}$ the first and second derivatives of $S$ with respect to a 
coordinate $\vp$ on $\partial M$. 

\begin{theorem} \label{thm: S close to const}
Let $(M,\bg)$ be a cuspidal surface of order $k$ with connected boundary. Assume that $S$ is a Morse function, that 
$S_{\vp\vp}< a_k$ on $\partial M$ and that $S_{\vp\vp}$ never takes the value $2$ at any local minimum, 
where $\vp$ is an arc length parameter on $\partial M$. Then
the exponential map \eqref{eqn:exp map general} is a homeomorphism for suitable $\tau_0$, $U$.
\end{theorem}
In particular, there is a neighborhood of $\partial M$ such that $U\setminus\dM$ is foliated by geodesics 
starting at $\dM$.

The exponential map  \eqref{eqn:exp map general} extends to the boundary by letting $\expdM(q,0)$ be the starting point of $\gamma_q$.
However, unlike in the case of constant $S$, the extension $\bS^1\times [0,\tau_0)\to U$ is not continuous 
in $q$ if $S$ is a Morse function; this is clear since the image of $\bS^1\times\{0\}$ is a finite set with at least two elements, a maximum and minimum of $S$.

The value $a_k$ in Theorem \ref{thm: S close to const} is optimal in the following sense.
\begin{theorem} \label{thm: S far from const}
Let $(M,\bg)$ be a cuspidal surface of order $k$ with connected boundary. Assume that $S$ is a Morse function and that 
$S_{\vp\vp} > a_k$ at some minimum of $S$, for an arc length parameter $\vp$ on $\dM$. Then the exponential 
map \eqref{eqn:exp map general} 
is not injective for any $\tau_0>0$.
\end{theorem}
That is, in any neighborhood $U$ of $\dM$ there are points through which at least two geodesics starting 
at the boundary pass.

Theorems \ref{thm: S close to const} and \ref{thm: S far from const} are illustrated in Figure \ref{fig:thms}, 
middle and right.

We also show that the conditions on $S$ in Theorem \ref{thm: S close to const} are satisfied if $M$ arises 
from a cuspidal surface singularity $X$ and if $\partial\Xtilde$ is contained in the boundary of a strictly convex subset 
of $\R^{n-1}$ which contains the origin, and is never doubly tangent to a sphere centered at the origin. See Theorem \ref{thm:convex bound}.

Theorems \ref{thm:S=const}, \ref{thm:surfaces}, \ref{thm: S close to const} and \ref{thm: S far from const}
may be summarized as follows, in the case of surfaces: There is a well-defined exponential map. For constant 
$S$ it is a smooth diffeomorphism near the boundary.  If $S$ is not too far from constant then it is a 
homeomorphism, though not on the boundary. If $S$ is far from constant then it may be not injective and also 
discontinuous.

\subsection*{Main ideas, outline of the proofs}
To simplify the notation, we assume in this outline that $k=2$. Thus the metric on a neighborhood $U$ of the boundary of $M$ is 
\begin{equation}
  \label{eqn:cusp metric inv}
\bg_{|U} = (1-2r^2 S + O(r^3))\,dr^2 + r^4 h\,.
\end{equation}
Let $\bg^*$ be the metric on the cotangent bundle $T^*\intM$ dual to $\bg$. Consider the energy function 
(Hamiltonian) $\Etilde=\frac12 \bg^*$ and the associated  Hamiltonian vector field $\Wtilde$ on $T^*\intM$. 
Then geodesics on $\intM$ are the projections to $\intM$ of integral curves of $\Wtilde$. Unit speed geodesics 
correspond to integral curves on the energy hypersurface $\{\Etilde=\frac12\}$. The degeneracy of $\bg$ at $\dM$ 
implies that $\Etilde$ and hence $\Wtilde$ are undefined over $\dM$. To make this explicit, let $m=\dim M$ and  $(r,\vp)$, 
$\vp=(\vp_1,\dots,\vp_{m-1})$ be local coordinates near a boundary point of $M$ and denote by 
$\xi,\eta=(\eta_1,\dots,\eta_{m-1})$ the dual coordinates on the fibers of $T^*M$. Since $\bg$ is a positive definite quadratic 
form in $dr$ and $r^2 \,d\vp$ whose coefficients are smooth functions of $r,\vp$ up to $r=0$, the function 
$\Etilde$ is a positive definite quadratic form in $\xi,\frac\eta{r^2}$ with coefficients smooth up to $r=0$. 
Therefore, rescaling
\begin{equation}\label{eqn:rescaling}
 \theta=\frac\eta{r^3},\quad E(r,\vp,\xi,\theta) = \Etilde(r,\vp,\xi,r^3\theta)
\end{equation}
yields a function $E$ which is smooth up to the boundary $r=0$, and simple calculations using the specific 
form \eqref{eqn:cusp metric inv} of $\bg$ show that the associated Hamiltonian vector field $\bW$ 
(which is $\Wtilde$ written in coordinates $r,\vp,\xi,\theta$) is $\frac1r$ times a vector field $\bV$ 
which is smooth up to the boundary $r=0$ and also tangent to the boundary.

Clearly, $\bV$ and $\bW$ have the same integral curves up to time reparametrization, so we need to analyze 
the integral curves of the rescaled vector field $\bV$.
It is essential for our analysis that  $\bV$ is sufficiently non-degenerate at its singular points  to make a precise analysis of its integral curves possible. For example,
the singular points of $\bV$ are hyperbolic if $S$ is a Morse function.

It may seem strange to use the rescaling \eqref{eqn:rescaling} instead of $\mu=\frac\eta{r^2}$, which is more 
naturally associated to cuspidal manifolds of order $2$ and is sufficient to make the energy a smooth function. But it turns out that the 
Hamilton vector field, written in coordinates $r,\vp,\xi,\mu$, needs to be multiplied by $r^2$ instead of $r$ 
to yield a smooth vector field, and that the resulting vector field is highly degenerate near its singular 
points, so a precise analysis via linearization would not be possible.

The rescaling \eqref{eqn:rescaling} may be given an invariant description as follows. This is the natural setting for our results and methods, but is not strictly necessary to 
understand most of the paper. Introducing $\theta$ corresponds to replacing the vector bundle $T^*M$ by a rescaled cotangent 
vector bundle which we denote by $\threeT^*M$, defined as the unique vector bundle over $M$  whose space of 
sections is the set of smooth one forms $\alpha$ on $M$ satisfying $\alpha(V)=O(r^3)$ for all smooth vector 
fields $V$ tangent to third order to the boundary, i.e. satisfying $dr(V)=O(r^3)$.\footnote{For simplicity, we 
assume a fixed boundary defining function $r$ is chosen, although to define $\threeT^*M$  it is enough to require 
the choice of $r$ up to changes of the form $r\mapsto a_1 r + a_2 r^2 + O(r^3)$ where $a_1>0$, $a_2$ are constants.} 
In coordinates, this space of sections is spanned by $dr$, $r^3d\vp_1,\dots,r^3 d\vp_{m-1}$ over $C^\infty(M)$.
Then $E$ is a smooth function on $\threeT^*M$ and the rescaled geodesic vector field $\bV$ is a smooth 
vector field on $\threeT^*M$.

We also need to describe the relation of $\threeT^*M$ to $T^*\partial M$. The rescaled cotangent bundle $\threeT^*M$ is a manifold with boundary $\partial(\threeT^*M) = \threeT^*_{\dM}M$, 
and $\pi^*r$ (where $\pi:\threeT^*M\to M$ is the projection), denoted $r$ for short, is a boundary defining function. 
The coordinate $\xi$ is invariantly defined on $\threeT^*_{\dM}M$, and for any $\xi_0\in\R$ the affine subbundle  
$\xi=\xi_0$  of $\threeT^*_{\dM}M$ may be naturally identified with 
$T^*\dM$. This identification is given in coordinates by $\xi_0\,dr + \theta\,r^3 d\vp \mapsto \theta d\vp$, 
where we write $\theta d\vp=\sum_{i=1}^{m-1}\theta_i d\vp_i$. That is, in coordinates the embedding 
$T^*\dM\to\threeT^*_{\dM}M$ is simply given by $(\vp,\theta) \mapsto (0,\vp,\xi_0,\theta)$ where on the 
left $\theta$ plays the role of the fiber coordinate on $T^*\dM$.

Our aim is to understand the dynamics of $\bV$ close to the boundary, so the dynamics of the restriction of $\bV$ 
to the boundary plays an important role in our analysis. 
Unit speed geodesics leaving the boundary have $\xi=1$ there, so this restriction may be considered as a 
vector field on $T^*\dM$ under the identification just described. We will see in Section \ref{sec:boundary dynamics} that the boundary dynamics is given by a 
damped Hamiltonian system with Hamiltonian
\begin{equation}
\label{eqn:Ebound intro} 
 \Ebound = S + \tfrac12 \gdM^*\quad \text{ on }T^*\dM .
\end{equation}
More precisely, in local coordinates it is given  by the equations
\begin{equation}
 \label{eqn:bd dynamics}
 \vpdot = \Ebound_\theta, \quad \thetadot = -\Ebound_\vp - 3\theta .
\end{equation}
Here $\Ebound_\theta=(\frac{\partial\Ebound}{\partial\theta_1},\dots,\frac{\partial\Ebound}{\partial\theta_{m-1}})$ etc.
The Hamiltonian $\Ebound$ corresponds to a particle moving on the Riemannian manifold $(\dM,\gdM)$ in the potential $S$.

The theorems now follow from a precise analysis of $\bV$. Singular points $\overbp\in \threeT^*_{\dM}M$ 
of $\bV$ (that is, points where $\bV$ vanishes) correspond to critical points $\bp\in\dM$ of $S$ via projection, 
and integral curves of $\bV$ leaving $\overbp$ correspond to geodesics starting at $\bp$. These integral 
curves foliate the unstable manifold $M^u_\bp\subset\threeT^*M$ of $\bp$. If $S$ is constant then each 
$\bp\in\dM$ is a critical point, and $\bV$ is transversally hyperbolic with respect to the submanifold of 
$\threeT^*M$ formed by the corresponding points $\overbp$, and the unstable manifold theorem yields  
Theorem \ref{thm:S=const}. For general $S$, the analysis of the exponential map breaks up into three parts:
\begin{enumerate}
 \item[a)] Understand the flow on $M^u_\bp$ for each critical point $\bp$.
 \item[b)] Understand  the regularity of $M^u = \bigcup_\bp M^u_\bp$, that is, how the various $M^u_\bp$ fit together.
 \item[c)] Determine whether $M^u$ projects diffeomorphically to $M$ in a neighborhood of $\dM$, under the 
projection $\threeT^*M\to M$.
\end{enumerate}
 If $S$ is a Morse function then the singular points of $\bV$ are isolated, and the linear part of $\bV$ 
is invertible at each singular point; this allows to do a) and answer c) affirmatively (for $M^u_\bp$) near 
each $\bp$. The main task now is to understand the global boundary dynamics, i.e.\ the global behavior of the 
unstable manifolds $\Gamma^u_\bp=\partial M^u_\bp$ of $\bV$ restricted to the boundary. In the case of surfaces 
this global analysis can be done.  This allows to do a), b) and c) in a full neighborhood of the boundary. Here 
some technicalities about linearizations enter; the assumption $S_{\vp\vp}\neq 2$ at  minima in Theorem \ref{thm: S close to const} is a non-resonance condition needed to ensure existence of a $C^1$ linearization. However, we believe that the theorem is true without this assumption. From this we then deduce Theorems \ref{thm: crit pts},
\ref{thm:surfaces}, \ref{thm: S close to const}, and \ref{thm: S far from const}.

For cuspidal metrics of order $k\geq2$ the only modifications to this outline are that $\tilde E$ is a positive definite quadratic form in $\xi, \frac{\eta}{r^k}$ and that $\theta = \frac\eta{r^{2k-1}}$, so that  the rescaled geodesic vector field $\bV$ lives naturally on the space
${}^{(2k-1)}T^*M$. In \eqref{eqn:Ebound intro}  one needs to replace $S$ by $\frac12k(k-1)S$ and in  \eqref{eqn:bd dynamics} one needs to replace $3\theta$ by $(2k-1)\theta$.

\subsection*{Structure of the paper} In Section \ref{sec:cusp metrics} we introduce and discuss cuspidal manifolds.
In Section \ref{sec:geodesic vector field} we calculate the rescaled geodesic vector field $\bV$ and analyze 
its singular points, and prove Theorem \ref{thm:S=const}. In Section \ref{sec:local behavior} we analyze the 
local behavior of the flow of $\bV$ for a Morse function $S$, and discuss linearizations. In Section
\ref{sec:boundary dynamics} we analyze the boundary dynamics, which we use in Section \ref{sec:exponential map} 
to prove the remaining theorems, after defining the exponential map. 
In Section \ref{sec:cusp sing} we apply our results to cuspidal singularities by first showing that the restriction of an ambient metric to a cuspidal singularity yields a cuspidal metric upon resolution, and then interpreting the quantities $S$ and $\bg_{\dM}$ in this context; we use this to prove Theorem \ref{thm:convex bound}.
In Section \ref{sec:examples} we give some 
examples. In the interest of better readability we present the arguments in Sections \ref{sec:geodesic vector field}-\ref{sec:exponential map} in detail for the case $k=2$ and then sketch the small adjustments needed for general $k$.

\subsection*{Further related work}
The exponential map  may be used to study differential geometric quantities near a singularity. 
For example, the description of the asymptotic behavior of the exponential map may be used to deduce the 
asymptotic behavior of the volume of balls centered at the singularity, as the radius tends to zero. In 
the case where $p$ is a smooth point, the coefficients in this expansion are related to the 
curvature at $p$. In \cite{Gri03} such an expansion was derived for real analytic isolated surface singularities; 
however, the balls were defined extrinsically, i.e. with distance 
defined as distance in the ambient space $\R^n$. 
A study of conical metrics related to our discussion of cuspidal metrics is done in \cite{Gri:NDOCS}.
The bi-Lipschitz geometry of (singular) algebraic subsets 
of $\R^n$ was investigated extensively  in \cite{Mos, Par, Bir, Va1, Va2}. However, very little is known about the 
precise intrinsic geometry of these sets, and the present paper provides a first insight into this question.
Global aspects of the geodesics of spaces with isolated singularities were studied in \cite{Ghi:GCM}.
\section{Cuspidal metrics} \label{sec:cusp metrics}

In this section we define cuspidal metrics on an $m$-dimensional manifold $M$ with boundary and show that they define invariantly a function $S$ and a metric $\bg_{\partial M}$ on the boundary.

\begin{definition}
Let $k\geq2$, $k\in\N$.  A {\em cuspidal metric of order $k$} on a manifold with compact boundary $M$ is a smooth semi-Riemannian metric $\bg$ on $M$ 
for which near any boundary point there are coordinates $(r,\varphi_1,\dots,\varphi_{m-1})$, and  the 
boundary defined by $r=0$, in which $\bg$ takes the form
\begin{equation}\label{eqn:def-cusp-metric}
\bg = \left(1-k(k-1)r^{2k-2} S + O(r^{2k-1})\right) \, dr^2 + r^{2k} \left( 2\sum_{i=1}^{m-1}b_i \,dr\,
d\varphi_i + \sum_{i,j=1}^{m-1} c_{ij}\,d\varphi_i\, d\varphi_j \right)
\end{equation}
where $S=S(\varphi)$, $b_i=b_i(r,\varphi)$ and $c_{ij}=c_{ij}(r,\varphi)$ are smooth and the bilinear form  $\sum_{i,j=1}^{m-1} c_{ij}(0,\varphi)\,d\varphi_i\, d\varphi_j$
is positive definite for each $\varphi$, so defines a Riemannian metric on the boundary.

A {\em cuspidal manifold} is a manifold with compact boundary endowed with a cuspidal metric.
\end{definition}
The main feature here are the powers of $r$ that occur. The reason for inserting the factor $k(k-1)$ is that in the case of cuspidal singularities the function $S$ has a natural meaning, see \eqref{eqn:S gdXtilde embedded}.
Observe that the form of metric  \eqref{eqn:def-cusp-metric} is unaffected by a change of the $\varphi$-coordinates. 
However, this is not true for the coordinate (boundary defining function) $r$. More precisely, we have:
\begin{lemma}\label{lem:cusp metric}
 A cuspidal metric of order $k$ on $M$ determines invariantly the function $S$ on $\partial M$, up to additive constants, 
as well as a Riemannian metric $\gdM$ on $\partial M$, given by $\gdM=\sum_{i,j=1}^{m-1} c_{ij}(0,\varphi)
\,d\varphi_i\, d\varphi_j$ in coordinates.
The metric also fixes the boundary defining function $r$ up to replacing it by $r+br^{2k-1}+O(r^{2k})$, where $b$ is 
constant on each connected component of $\dM$. 
\end{lemma}
\begin{proof}
Let $R$ be another boundary defining function. We find conditions on $R$ ensuring that the metric has the form \eqref{eqn:def-cusp-metric} when written in terms of $R$ instead of $r$. Write 
$r=w_1(\varphi)R + \dots + w_{2k-1}(\varphi) R^{2k-1}  +O(R^{2k})$.
Then
$$ dr = \sum_{i=1}^{2k-1} i w_i R^{i-1}\, dR + \sum_{i=1}^{2k-1} R^i dw_i + O(R^{2k-1})\,dR + O(R^{2k}) d\varphi\,.$$
% $$dr=\left(w_1+\dots+(2k-1)w_{2k-1} R^{2k-2} + O(R^{2k-1})\right)dR + \left(R dw_1 + \dots + R^{2k-1}dw_{2k-1}\right) + O(R^{2k}).$$
 Inserting this in \eqref{eqn:def-cusp-metric} we see that the coefficient of $dR^2$ in $\bg$ will be of the form
$1+O(R^{2k-2})$ iff $w_1\equiv 1$, $w_2\equiv\dots\equiv w_{2k-2}\equiv 0$, and then 
%$dr=(1+3R^2 w_3)\,dR + R^3 dw_3$ and therefore
%\begin{equation}
\begin{equation}\label{eqn:metric change}
\begin{gathered}
 \left(1-k(k-1)r^{2k-2} S + O(r^{2k-1})\right)\,dr^2 
 = A\,dR^2 + 2R^{2k-1}dR \,dw_{2k-1} + O(R^{4k-3}) \\
A = 1 + R^{2k-2}\left(-k(k-1) S + 2(2k-1)w_{2k-1}\right) + O(R^{2k-1}) \,.
\end{gathered}
\end{equation}
%\end{equation}
This shows that the mixed term involving $dR d\varphi$ in $\bg$ is $O(R^{2k})$ (as required in \eqref{eqn:def-cusp-metric}) iff $dw_{2k-1}=0$, i.e.\ $w_{2k-1}$ is locally constant.
Therefore, under a change of the boundary defining function $r$ preserving the form \eqref{eqn:def-cusp-metric} of the metric, the coefficient of $r^{2k-2}\,dr^2$ (hence $S$) is determined up to a locally constant term on $\partial M$ and $r$ is determined to the order claimed. Also, the $d\varphi_id\varphi_j$ term is changed only by $O(R^{4k-2})$ and transforms like a Riemannian metric under changes of the $\varphi$-coordinates, so the metric on $\partial M$ is well-defined.
\end{proof}
%%%\ownremark{Question to be pondered sometime: How can $S$, or rather $dS$, be defined without coordinates?}

%
%
%
%
%
%
%
%
%
%
%
%
%
%
%
%
%
\section{The rescaled geodesic vector field}
\label{sec:geodesic vector field}
In this section we analyze the rescaled geodesic vector field $\bV$ for a cuspidal metric \eqref{eqn:def-cusp-metric}. We show that $\bV$  is smooth, find its singular locus and its linearization at any singular point.
 This is the foundation for the proofs of all theorems. At the end of the section we prove Theorem \ref{thm:S=const}.

We present the detailed computations for $k=2$ and indicate the changes required for general $k$ at the end of the section. 
We always work in coordinates $(r,\vp)$ in which \eqref{eqn:def-cusp-metric} holds and write the metric in simplified notation as
$$ \bg = (1-2r^2S + O(r^3))\, dr^2 + r^4 2B\,dr\,d\varphi + r^4C\,d\varphi^2 $$
where $d\vp=(d\vp_1,\dots,d\vp_{m-1})$, $B=(b_1,\dots,b_{m-1}),C=(c_{ij})_{i,j=1,\dots,m-1}$ (with $c_{ij}=c_{ji}$) are smooth functions of $r,\varphi$ 
and $S$ is a smooth function of $\varphi$. All vectors will be treated as column vectors. By assumption $C>0$ at 
$r=0$. Here and in the sequel, $O(r)$ denotes $r$ times a smooth, possibly vector- or matrix-valued function defined near $\dM$, hence can be differentiated.

Let $\xi$, $\eta$ denote the cotangent coordinates dual to $r$, $\vp$, respectively. In order to calculate 
the metric dual to $\bg$, we need to invert the coefficient matrix of $\bg$. We use the general formula 
(for a scalar $a$, vector $b$ and symmetric matrix $c$)
$$ 
\begin{pmatrix}
 a & b^t \\
 b & c
\end{pmatrix}
^{-1}
= d^{-1}
\begin{pmatrix}
 1 & -b^t c^{-1} \\
 -c^{-1}b & c^{-1}d + c^{-1}bb^tc^{-1}
\end{pmatrix}
, \quad
d = a - b^tc^{-1}b
$$
which can be verified by direct calculation, with $a=1-2r^2S+O(r^3)$, $b=r^4B$, $c=r^4C$. Then $d=1-2r^2S + O(r^3)$, 
hence $d^{-1}=1+2r^2 S + O(r^3)$, and we get
\begin{align*}
 \bg^* &= \left(1+2r^2S + O(r^3)\right) \xi^2 + 2\left(- B^tC^{-1} + O(r^2)\right)\xi\eta + r^{-4} \left(C^{-1} + O(r^4)\right) \eta^2 \\
 &= \xi^2 + 2r^2 G(r,\vp,\xi,\frac\eta{r^3})
\end{align*}
 ($C^{-1} \eta^2$ is to be interpreted as $\eta^t C^{-1} \eta$)
where 
\begin{equation}
 \label{eqn:def G}
\begin{gathered}
 G(r,\vp,\xi,\theta) = G_0(r,\vp,\xi,\theta) + 
  \Gtilde(r,\vp,\xi,\theta),\\ G_0 = S(\vp)\xi^2 + \tfrac12 \Cpartial(\vp)^{-1}\theta^2, \qquad
    \Gtilde = O(r)\xi^2 + O(r)\xi\theta + O(r)\theta^2
 \end{gathered}
\end{equation}
where $\Cpartial(\vp) = C(0,\vp)$.
% \ownremark{note the change from a previous version: in $g^{-1}$ the error term in the $\eta^2$ coefficient 
% is $r^4$, not $r^2$, resulting in $r^4$ in the previous formula; also note that $O(r)\xi\theta$ could be
% replace by $O(r^2)\xi\theta$, see own remarks below}
The energy $\Etilde(r,\vp,\xi,\eta)$ is defined as half the dual metric, so if we rewrite it in 
the rescaled coordinates \eqref{eqn:rescaling}
 we get
\begin{equation}
\label{eq: E}
E(r,\vp,\xi,\theta) = \frac12\xi^2 + r^2G(\vp,r,\xi,\theta)
\end{equation}
The geodesics are projections of integral curves of the Hamiltonian vector field associated with the
energy function $E$ and the symplectic form $dr\wedge d\xi + d\vp\wedge d\eta$, that is they are 
the $r,\vp$ part of the solutions of the following differential equation:
\begin{equation}
\label{eq: diffeq E}
\begin{aligned}
r' &= \Etilde_\xi  & \xi' &= - \Etilde_r \\
\vp' &= \Etilde_\eta & \eta' &= - \Etilde_\vp
\end{aligned}
\end{equation}
In order to rewrite this using the $\theta$ coordinate, we differentiate \eqref{eqn:rescaling} with 
respect to the various variables and obtain
\begin{equation*}
\begin{aligned}
E_\xi & = \Etilde_\xi &E_r& = \Etilde_r + 3r^2\theta \Etilde_\eta \\
 E_\theta & = r^3 \Etilde_\eta & E_\vp & = \Etilde_\vp \\
\end{aligned}
\end{equation*}
and $$ \theta^\prime  = \displaystyle
{\frac{\eta^\prime}{r^3} - 3 \frac{r^\prime}{r}\theta} $$
Therefore, in the new coordinates the differential equations \eqref{eq: diffeq E} are
\begin{equation}
\vspace{4pt}
\label{eq: diffeq E1}
\begin{aligned}
r' &= E_\xi  & \xi' &= - E_r +3\frac{\theta}{r} E_\theta \\
{\vp}' &= \frac{1}{r^3}E_\theta & {\theta}' &= - \frac{1}{r^3}E_\vp - 3\frac{\theta}{r}E_\xi
\end{aligned}
\vspace{4pt}
\end{equation}
This new system corresponds to a vector field $\bW$ that we call
\textit{the geodesic vector field},
which is well defined and smooth outside $\{r=0\}$. \\
Since $E = \frac12\xi^2 + r^2G$ with $G$ smooth as a function of  $r,\vp,\xi,\theta$ by \eqref{eq: E}, we obtain from \eqref{eq: diffeq E1} the following fundamental fact.
\begin{lemma}
 \label{lem:V smooth}
The rescaled geodesic vector field  $\bV=r\bW$ is smooth up to $r=0$, that is, it is a smooth vector field on the space $\threeT^*M$. In coordinates $r,\vp,\xi,\theta$ its integral curves are the 
solutions of the system
\begin{equation}
\vspace{4pt}
\label{eq: diffeq E2}
\begin{aligned}
\dot{r} &= r[\xi +  r^2G_\xi]  & \dot{\xi} &=  r^2[-2G - rG_r+ 3\theta G_\theta]\\
\dot{\vp} &=  G_\theta & \dot{\theta} &= -\left[ G_\vp + 3\xi\theta + 3r^2\theta G_\xi \right]
\end{aligned}
\end{equation}
where $G$ is defined in \eqref{eqn:def G}.
\end{lemma}

(Throughout the paper the time variables for integral curves of $\bW$, $\bV$ are denoted $\tau$, $t$, 
respectively, and $\frac d{d\tau}$ by a prime and $\frac d{dt}$ by a dot.)
Unit speed geodesics lie on the energy level $\{E=\tfrac12\}$. This is a $(2m-1)$-dimensional manifold which 
intersects the boundary $\{r=0\}$ transversally and whose intersection with the boundary has two connected 
components corresponding to $\xi=\pm 1$, by \eqref{eq: E}. In a neighborhood of this intersection each 
component can be parametrized by $(r,\vp,\theta)$. The $\rdot$ equation shows that $\xi>0$ for geodesics 
leaving the boundary. We will always consider this component. 

It is important to understand the behavior of the energy hypersurface $\{E=\tfrac12\}$ near $r=0$, and also the leading terms
of  $\bV$
near $r=0$. Note that the energy surface is non-compact since $\theta$ is not restricted 
to a bounded set (as opposed to $\eta$ or $\frac\eta{r^2}$).
For the following estimates it is important that they are uniform in all variables in a neighborhood of $\{ r=0\}$. 
In particular, $O(r)$ means $r$ times a smooth {\em bounded} function on $\{E=\tfrac12\}$. This makes the error estimates in the sequel a little cumbersome. It may be useful to read them separately for bounded $\theta$ -- which is all that is needed for analysis near the singular locus of $\bV$ -- and general $\theta$, which is needed in global arguments.

Since $\bg^*$ is a quadratic form in $\xi$ and $\frac{\eta}{r^2}=r\theta$ which is uniformly positive definite, 
\begin{equation}
 \label{eqn:xi rtheta bounded}
\xi,\ r\theta \text{ are bounded on } \{E=\tfrac12\}.
\end{equation}
Then \eqref{eqn:def G}, \eqref{eq: E} imply that
\begin{equation}
 \label{eqn: energy shell}
 \xi^2 + r^2\Cpartial^{-1}\theta^2 = 1 + O(r^2)+O(r^3)\theta^2 \qquad\text{ on } \{E=\tfrac12\}.
\end{equation}
% \ownremark{The $O(r^2)$ stems from the term $S\xi^2$ in $G_0$ and from the $O(r)\xi\theta$ term in 
% $\Gtilde$, which is of the same order as $\xi^2+r^2\theta^2$; as explained in the next 'own remark', 
% this term may be made to be $O(r^2)\xi\theta$ by a good choice of $\vp$ coordinate -- that is, by choosing 
% a good trivialization of $M$ near the boundary to first order --, and then \eqref{eqn: energy shell} may be 
% improved to $(1+2r^2S)\xi^2 + r^2 C^{-1}\theta^2 = 1 + O(r^3)$; however, there seems to be no use for this 
% improvement } 
Also $G_\xi=O(1)$, $\Gtilde_\theta=O(r)+O(r)\theta$, 
$|rG_r|+|\Gtilde_\vp| = O(r)+O(r)\theta^2$
%-- REMARK: maybe use chinese brackets $\langle\theta\rangle = 1+|\theta|$ here, then $\Gtilde_\theta = O(r\langle\theta\rangle), 
%|rG_r|+|\Gtilde_\vp| = O(r\langle\theta\rangle^2)$ --
(the two error terms are not 
comparable since $\theta$ is unbounded) and $3\theta G_\theta - 2G = 2C^{-1}\theta^2 + O(1) + O(r)\theta^2 = r^{-2}\left[2(1-\xi^2) 
+ O(r^2) + O(r^3)\theta^2\right]$ (using \eqref{eqn: energy shell}), so \eqref{eq: diffeq E2} yields, with uniform 
estimates on $\{E=\tfrac12\}$, 
\begin{equation}
\label{eq: diffeq E1 leading terms}
\begin{aligned}
\dot{r} &= r\xi +  O(r^3)  & \dot{\xi} &=  2(1-\xi^2) + O(r^2) + O(r^3)\theta^2\\
\dot{\vp} &=  (G_0)_\theta + O(r) +O(r)\theta & \dot{\theta} &= -(G_0)_\vp - 3\xi\theta +O(r) + O(r)\theta^2
\end{aligned}
\end{equation}
% \ownremark{The error term in the $\thetadot$ equation may be improved to $O(r)$, which may be useful for 
% some purposes, as follows: The only term which yields $O(r\theta)$ instead of $O(r)$ is $\Gtilde_\vp$ and 
% in this the $O(r)\xi\theta$ part on $\Gtilde$. But note that this will be $O(r^2)\xi\theta$ if $B=O(r)$, 
% i.e. $B=0$ at the boundary. Now if $\bg$ is any cuspidal metric, then one can easily arrange this. Simply 
% change the $\vp$ variable to $\psi$ where $\vp = \psi - r BC^{-1}$, for then $d\vp = d\psi - BC^{-1}dr + O(r)$ 
% and hence $r^4 2B drd\vp + r^4 Cd\vp^2 = r^4 2B dr d\psi + r^4 C(d\psi^2 -2BC^{-1}dr d\psi) + O(r^4)dr^2 + O(r^5)$, 
% and the $r^4$ part of the $drd\psi$ term cancels -- this is written for surfaces but should work the same
%  way in any dimension; note that then the $\Gtilde$ part has $\theta$ only appearing in the combination $r^2\theta$,
% which is $r\mu$ -- maybe useful for nearby geodesics}
\begin{lemma} \label{lem:sing points}
In a neighborhood of $\{r=0\}$ the singular locus of $\bV$ on the energy hypersurface $\{E=\tfrac12\}$ is 
$\{r=0,\theta=S_\vp =0\}$. 
\end{lemma}
\begin{proof}
By \eqref{eq: diffeq E1 leading terms}  we have $\xidot\geq 1-2\xi^2$ in a neighborhood of $\{r=0\}$, 
so $\xidot=0$ implies $|\xi|\geq \frac1{\sqrt2}$ and hence $|\dot{r}| \geq \frac1{\sqrt2} r+O(r^3)$, 
so at a singular point of $\bV$ we must have $r= 0$ and hence $\xi=\pm 1$. Then $\dot{\vp} = C^{-1}\theta = 0$ gives 
$\theta =0$, and then $\dot{\theta}=0$ gives $S_\vp =0$. 
\end{proof}
Therefore singularities of $\bV$ correspond to critical points of $S$, that is, 
to points $\bp=\bp(\vp_0)\in\partial M$  where $S_\vp(\vp_0)=0$. We denote the corresponding 
singular point of $\bV$ having $\xi=1$ and $\theta=0$ by 
$\overbp  \in \{E=\tfrac12\}$.
\begin{lemma}\label{lem:linearization}
 The linear part of $\bV$ on the energy hypersurface at a singular point $\overbp=(0,\vp_0,1,0)$ is
$$ \rdot = r,\quad \vpdot=\gamma_1 r + \Cpartial^{-1}\theta,\quad
\thetadot= \gamma_2 r - S_{\vp\vp} \cdot(\vp-\vp_0) - 3\theta $$
with $S_{\vp\vp},\Cpartial$ evaluated at $\bp$ and for suitable $\gamma_1,\gamma_2$.
\end{lemma}
The precise values of $\gamma_1,\gamma_2\in\R^{m-1}$ are inessential for the analysis. 
$\gamma_1=-B^tC^{-1}$
%%%\ownremark{which is zero in the improved coordinates from the previous ownremark}, 
and $\gamma_2$ is determined by $B,C$ and the $O(r^3)$ term in the coefficient of $dr^2$ in $\bg$.
We rewrite the linear part in matrix form:
\begin{equation}\label{eqn:linearization}
\text{Linear part of $\bV$:}\quad
\begin{pmatrix}
 \rdot \\ \vpdot \\ \thetadot
\end{pmatrix}
=
\begin{pmatrix}
 1 & 0 & 0 \\
 \gamma_1 & 0 & \Cpartial^{-1} \\
 \gamma_2 & -S_{\vp\vp} & -3\,\Id
\end{pmatrix}
\begin{pmatrix}
 r \\ \vp-\vp_0 \\ \theta
\end{pmatrix}
\end{equation}

We can now settle the case of constant $S$. 
\subsection*{Proof of Theorem \ref{thm:S=const}}
The proof is analogous to the proof of  the result of Melrose \& Wunsch  in the conical case \cite{MelWun}.
The idea is that by Lemma \ref{lem:sing points} for constant $S$ every point on the set $\calC=\{(r,\vp,\xi,\theta):\,r=0,\theta=0,\xi=1\}$ is a
singular point of $\bV$. The linearization in Lemma \ref{lem:linearization} shows that $\bV$ is transversally hyperbolic along $\calC$, so there is a unique integral curve starting at each point of $\calC$. Since $\calC$ projects diffeomorphically to $\partial M$, the projections of these curves give the desired foliation of a neighborhood of $\partial M$.

 Let $\overbp=(0,\vp_0,1,0)\in \calC$. The set $\{\vp=\vp_0\}$ is transversal to 
$\calC$, and projecting  the linear part of $\bV$ at $\overbp$ onto the tangent space of this set  we obtain
\begin{equation}
\label{eqn:linearization r theta} 
 \begin{pmatrix}
 1 & 0 \\
 \gamma_2 & -3\,\Id
\end{pmatrix}
\end{equation}
(delete the second block of rows and columns in the matrix \eqref{eqn:linearization}).
This has a single eigenvalue $1$ and the $(m-1)$-fold eigenvalue $-3$, so $\bV$ is normally hyperbolic 
along $\calC$. The unstable eigenvector $v$ (for the eigenvalue $1$) has non-vanishing $\partial_r$-component. 
By the stable manifold theorem \cite{HPS}
there exists a unique smooth unstable submanifold $N$ of dimension $m$ with $\partial N=N\cap\{r=0\}=\calC$, 
transversal to $\{r=0\}$, which is the union of trajectories of $\bV$ starting at points of $\calC$.
The restriction of $\bV$ to $N$ is a smooth vector field on $N$ which vanishes at $r=0$, hence 
$\bW=\frac1r\bV_{|N}$ extends smoothly from $N\cap\{r>0\}$ to all of $N$ by Taylor's theorem. Also, 
$\bW$ has a non-zero $\partial_r$-component at $\partial N$ by the $\rdot$ equation in \eqref{eq: diffeq E1}. 
Near $r=0$, $N$ is a graph over $(r,\vp)$ since its tangent space is spanned by $T\calC$ and $v$. Therefore, 
the projection $(r,\vp,\theta)\mapsto(r,\vp)$ restricts, near $\partial N$, to a diffeomorphism 
$\psi : N \to M$. 
Therefore, the vector field $\bW'=\psi_*(\bW|_N)$ is well-defined near $\partial M$, and its integral 
curves are geodesics in $M$. Since $\bW'$ is transversal to $\partial M$, the result follows 
from standard facts about flows.

\begin{remark}
1) This can be generalized to the case when $\Sbound$ is constant only on an open subset of $\dM$, with the same proof.
\\
2) By the standard Gauss lemma (which is also applicable here, cf. the footnote in the proof of the conical 
theorem in \cite{MelWun}), the theorem implies that for constant $S$ the normal form 
\eqref{eqn:def-cusp-metric} may be improved: There is a trivialization of a neighborhood $U$ 
of the boundary, $ \partial M\times [0,\tau_0)\cong U$ (given by the flow of the vector field $\bW'$ 
in the proof), so that 
$$ \bg_{|U} = dr^2 + r^4 h(r)$$
(resp.\ $dr^2+r^{2k}h(r)$ for cusps of order $k$)
where $h$ is a smooth family of metrics on $\partial M$.
\end{remark}  
\subsection*{Adjustments for general order $k$}
For cuspidal manifolds of order $k\geq2$ the correct rescaling of the cotangent variable is
$$ \theta = \frac\eta{r^{2k-1}} \,.$$
The same line of argument as above yields, generalizing  \eqref{eq: E}, \eqref{eqn:def G}, the energy 
$$E=\frac12\xi^2 + r^{2k-2}G\quad\text{where}\quad G=G_0+\Gtilde,\ G_0 = \tfrac12 k(k-1) S \xi^2 + \tfrac12 \Cpartial^{-1}\theta^2 $$
and $    \Gtilde = O(r)\xi^2 + O(r)\xi\theta + O(r)\theta^2$.
This yields the rescaled geodesic vector field, $\bV = r \bW$,
\begin{equation}
\vspace{4pt}
\tag{$\ref{eq: diffeq E2}_k$}
\label{eq: diffeq E2 k}
\begin{aligned}
\dot{r} &= r[\xi +  r^{2k-2}G_\xi]  & \dot{\xi} &=  r^{2k-2}[-(2k-2)G - rG_r+ (2k-1)\theta G_\theta]\\
\dot{\vp} &=  G_\theta & \dot{\theta} &= -\left[ G_\vp + (2k-1)\xi\theta + (2k-1)r^{2k-2}\theta G_\xi \right]
\end{aligned}
\end{equation}

Then $\xi$ and $\frac\eta{r^k}=r^{k-1}\theta$ are bounded on the energy shell, and \eqref{eqn: energy shell} generalizes to
\begin{equation}
 \label{eqn: energy shell k}
 \tag{$\ref{eqn: energy shell}_k$}
 \xi^2 + r^{2k-2}\Cpartial^{-1}\theta^2 = 1 + O(r^{2k-2})+ O(r^{2k-1})\theta^2 \quad\text{ on }\{E=\tfrac12\}\,.
\end{equation}
Then \eqref{eq: diffeq E2 k} yields the equations near the boundary, 
generalizing \eqref{eq: diffeq E1 leading terms},
\begin{equation}\label{eq: diffeq E1 leading terms k}
\tag{$\ref{eq: diffeq E1 leading terms}_k$}
\begin{aligned}
\dot{r} &= r\xi +  O(r^{2k-1}) + O(r^{2k})\theta  & \dot{\xi} &=  k(1-\xi^2) + O(r^{2k-2})+O(r^{2k-1})\theta^2 \\
\dot{\vp} &=  (G_0)_\theta + O(r) + O(r)\theta & \dot{\theta} &= -(G_0)_\vp - (2k-1)\xi\theta + O(r) + O(r)\theta^2
\end{aligned}
\end{equation}
At a singular point this has linear part on the energy surface as in \eqref{eqn:linearization} resp.\ \eqref{eqn:linearization r theta}
with 
\begin{equation}
 \label{eqn:replace k}
\text{$S$ replaced by $\tfrac12 k(k-1) S$ and $3\Id$ by $(2k-1)\Id$.}
\end{equation}
The proofs of the lemmas and of Theorem \ref{thm:S=const} proceed as in the case $k=2$.
\section{Local behavior of the rescaled geodesic vector field near singularities}
\label{sec:local behavior}
Let $\overbp$ be a singular point of the rescaled geodesic vector field $\bV$ corresponding to a 
critical point $\bp$ of $\Sbound$. 
We will determine  the local phase portrait of  $\bV$ near $\overbp$ in the case 
where $\bp$ is a non-degenerate critical point of $S$. 

We first determine the eigenvalues and eigenvectors of the linear part  \eqref{eqn:linearization} of 
$\bV$ at $\overbp$. Recall that the Hessian of $S$ at a critical point $\bp$ is a well-defined symmetric 
bilinear form on $T_\bp\partial M$. By way of the metric $\gdM$, this form corresponds to a linear 
endomorphism of $T_\bp\partial M$, which we call the {\em Hessian of $S$ relative to $\gdM$} at $\bp$.
%%\notevg{I put the lines back in, for the purpose of explanation to people less comfortable with the previous statement}
For example, if one chooses coordinates in which $\partial_{\vp_1},\dots,\partial_{\vp_{m-1}}$ 
are orthonormal with respect to $\gdM$ at $\bp$ then this is simply the matrix of 
partial derivatives $S_{\vp_i \vp_j}$, and $\Cpartial=\Id$ at $\bp$. 
A simple calculation using \eqref{eqn:linearization} (considering also \eqref{eqn:replace k} for general order $k$) yields:
\begin{lemma}\label{lemma-eigenvalues}
The linear part of $\bV$ at a singular point $\overbp$ has an eigenvalue $\lambda_1(\bp)= 1$ with 
eigenvector transverse to $r=0$, and for each eigenvalue $a$ of the Hessian of 
$S$ relative to $\gdM$ at $\bp$ and each corresponding eigenvector 
$u\frac{\partial}{\partial \vp}=\sum_{i=1}^{m-1}u_i\frac{\partial}{\partial \vp_i}$ the two eigenvalues (in the case $k=2$)
\begin{equation}
 \label{eqn:eigenvalues}
\displaystyle{\lambda_2 (\bp,a)= \frac{-3+\sqrt{9-4a}}{2}}\ \text{ and }\ 
\displaystyle{\lambda_3 (\bp,a)= \frac{-3-\sqrt{9-4a}}{2}} 
\end{equation}
with eigenvectors
\begin{equation}
 \label{eqn:eigenvectors}\displaystyle{\nu_2 (\bp,a,u)= u \frac{\partial}{\partial \vp} +\lambda_2(\bp,a) (\Cpartial u)
\frac{\partial}{\partial \theta}} \text{ and } \displaystyle{\nu_3 (\bp,a,u)= u \frac{\partial}{\partial \vp} +
\lambda_3(\bp,a)(\Cpartial u)\frac{\partial}{\partial \theta}}
\end{equation}
For general $k$ the eigenvalues are
$$ \lambda_{2/3}(\bp,a)=\frac{-(2k-1) \pm \sqrt{(2k-1)^2 - 2k(k-1) a}}2$$
with the same eigenvectors as before.
\end{lemma}
Note that the eigenvalues $\lambda_{2,3}(\bp,a)$ are non-real if $a > a_k$ where
$$ a_k = \frac{(2k-1)^2}{2 k(k-1)} = 2 + \frac 1{2k(k-1)}.$$
In particular $a_2=9/4$.

This allows us to analyze the geodesics near non-degenerate critical points. We only note the cases 
of maxima and minima of $S$, other critical points can be treated similarly.
\begin{proposition}\label{prop-localmax}
Let $\bp$ be a non-degenerate critical point of $\Sbound$ on a cuspidal manifold $M$.
\begin{enumerate}
\item
If $\Sbound$ has a local maximum at $\bp$ then there is a neighborhood $U$ of $\bp$ in $M$ so 
that $U\cap\intM$ is foliated by geodesics.
\item
If $\Sbound$ has a local minimum at $\bp$ then there is a unique geodesic which starts at $\bp$.
\end{enumerate}
\end{proposition}
\begin{proof}
At a local maximum of $\Sbound$, all eigenvalues $a$ of its Hessian are negative, so all eigenvalues 
of the linear part of $\bV$ 
are real and $\lambda_3(\bp,a) < 0 < \lambda_2(\bp,a)$. By the unstable manifold theorem, the unstable manifold 
$M_\bp^u$ is a smooth $m$-dimensional submanifold of $\threeT^*M$ whose tangent space at $\overbp$ is 
spanned by the eigenspaces 
corresponding to $\lambda_1$ and the various $\lambda_2(\bp,a)$. Eigenvectors for $\lambda_1$ have non-zero 
$\partial_r$-component, and the $\vp$-components of the eigenvectors $\nu_2(\bp,a)$ span $T_\bp\partial M$. 
Therefore the projection of
$M_\bp^u$ onto $M$ is a diffeomorphism in a neighborhood of
$\bp$. Since $M_\bp^u\setminus \overbp$ is locally foliated by trajectories of $\bV$, the result for a 
maximum follows.

At a minimum of $\Sbound$, all values $a$ are positive, so the linear part of $\bV$ has a unique eigenvalue, 
$\lambda_1=1$, with positive real part. By the unstable manifold theorem, there is a unique trajectory 
leaving $\overbp$. The result follows.
\end{proof}

\medskip

From now on we focus on cuspidal surfaces. In this case the Hessian of $S$ at $\bp$ has only one eigenvalue $a$, which we denote by 
$$ a(\bp) = S_{\vp\vp}(\bp) \Cpartial (\bp)^{-1}$$
In arc length coordinates this is simply $S_{\vp\vp}(\bp)$. We write $\lambda_{2,3}(\bp) = 
\lambda_{2,3}(\bp,a(\bp))$ and $\nu_{2,3}(\bp) = \nu_{2,3} (\bp,a(\bp),1)$.
See \eqref{eqn:meaning Spp} for the geometric meaning of $S_{\vp\vp}$ and hence $a$ in the case of cusp singularities.

\subsection*{Linearization near a singular point}
A triple of 
%pairwise distinct 
complex numbers $\lambda_1,\lambda_2,\lambda_3$ are \em resonant \em
if there exist $i\in \{1,2,3\}$ and $a_1,a_2,a_3 \in \N_0$ 
with $a_1 + a_2 + a_3 \geq 2$ such that 
\vspace{4pt}
\begin{equation}\label{eq:resonances}
\vspace{4pt}
a_1\lambda_1 + a_2\lambda_2 + a_3\lambda_3  = \lambda_i.
\end{equation}
This relation, or the triple $(a_1,a_2,a_3)$, is called \em a resonance relation of type $i$\em\  
among $\lambda_1$, $\lambda_2$ and $\lambda_3$.

\smallskip
Given  a non degenerate local extremum $\bp$ of the function $S$, we observe that we always have the relation 
$(2k-1)\lambda_1 (\bp)+ \lambda_2 (\bp) + \lambda_3(\bp) = 0$. Multiplying by any natural number and adding 
$\lambda_i(\bp)$ we obtain three families of resonance relations.

We recall the following very useful linearization criterion (Samovol's criterion).

\begin{theorem}[\cite{Sam}]\label{TheoremLinearizationSamavol}
Let $(\mu_1 \ldots \mu_s,\sigma_1,\ldots,\sigma_u)$ be the set of eigenvalues of the
linear part of a smooth vector field $V$ at $\bo$, a singular point of $V$. Assume that
\begin{center}
\vspace{4pt}
$Re(\mu_s)\leq \ldots \leq Re (\mu_1) < 0 < Re(\sigma_1) \leq \ldots \leq Re(\sigma_u)$
\vspace{4pt}
\end{center}
Assume there exists a positive integer $l$ such that for each resonant relation,
\begin{center}
\vspace{4pt}
$\lambda = \sum r_i^- \mu_i + \sum r_j^+\sigma_j$,
\vspace{4pt}
\end{center}
where $\lambda$ is an eigenvalue, there exists an integer $i \leq s$ or $j \leq u$ such that
\begin{center}
\vspace{4pt}
$- l Re(\mu_i)  < - Re(r_1^- \mu_1 + \ldots + r_i^-\mu_i)$ or
$l Re(\sigma_j)  < Re(r_1^+ \sigma_1 + \ldots + r_i^+\sigma_j)$.
\vspace{4pt}
\end{center}
Then the vector fields $V$ is $C^l$-conjugate to its linear part in a neighborhood of $\bo$.
\end{theorem}

As a consequence of this result and of the remark, we get the following proposition.
\begin{proposition}\label{linearization-loc-min}
Assume that $M$ is a cuspidal surface of order $k$ and that $\bp\in \partial M$ is a non-degenerate local minimum 
for the function $\Sbound$ for which $a(\bp)<a_k$ and $a(\bp)\neq 2$. If $\lambda_2(\bp)\in\Q$ 
then $\bV$ is $C^1$ linearizable near $\overbp$, otherwise it is $C^{2k-2}$ linearizable.
\end{proposition}
\begin{proof}
%We give the details for $k=2$, the general case is similar. 
Note that
  $a(\bp) \neq 2$ is equivalent to $\lambda_2(\bp) \neq -(k-1)$.

Write $\lambda_i=\lambda_i(\bp)$.
The assumptions imply $\lambda_1=1$, $\lambda_3<-\frac{2k-1}{2}<\lambda_2<0$ and $\lambda_2 \neq -(k-1)$.

1) Assume $\lambda_2$ is irrational. The resonance relations between the eigenvalues 
are of the form $\delta (2k-1,1,1) + R$, where $R \in \{(1,0,0),(0,1,0),(0,0,1)\}$ and $\delta
\in \N$. For each such relation we always find
\begin{center}
\vspace{4pt}
$(2k-2)\lambda_1 < [(2k-1)\delta +\omega]\lambda_1$ for $\omega = 0$ or $1$.
\vspace{4pt}
\end{center}
We conclude using Theorem \ref{TheoremLinearizationSamavol} that $\bV$ is $C^{2k-2}$ linearizable at $\overbp$.

\medskip 
2) Assume that $\lambda_2 = -\frac{p}{q} \neq -(k-1)$ is a negative rational number written in its irreducible form.

\smallskip
Let $(a_1,a_2,a_3)$ be a resonance relation of the form $a_1 \lambda_1 +a_2\lambda_2 + a_3\lambda_3 = \lambda_1$.
\\
Then $a_1>1$ since otherwise $0>a_2\lambda_2 + a_3\lambda_3 =(1-a_1)\lambda_1 \geq 0$. Hence 
$a_1 \lambda_1 > \lambda_1$. 

So for all resonance relations of type 1 we always find $a_1 \lambda_1 > \lambda_1$, so that we can take 
$j=1,l=1$ in regards of Samovol's criterion. 

\medskip
Let $(a_1,a_2,a_3)$ be a resonance relation of the form $a_1 \lambda_1 +a_2\lambda_2 + a_3\lambda_3 = \lambda_2$.
\\
If $a_2 \geq 2$ then $a_2(-\lambda_2)> -\lambda_2$. 
\\
If $a_2 = 1$ then $a_1 + a_3 \lambda_3 = 0$, thus $a_1 \geq k$ and so $a_1 \lambda_1 > \lambda_1$.
\\
If $a_2 = 0$ we check that $(a_1,a_3) \neq (1,1)$ since $\lambda_2 \neq -(k-1)$. Thus we either have $a_1 \lambda_1 > \lambda_1$
or $a_2(-\lambda_2)+ a_3 (-\lambda_3) > -\lambda_3$.
   
So for all resonance relations of type 2 we can take $l=1$ and, respectively, $i=1$, $j=1$, $j=1$ or $i=2$, 
in view of Samovol's criterion. 

\medskip
Let $(a_1,a_2,a_3)$ be a resonance relation of the form $a_1 \lambda_1 +a_2\lambda_2 + a_3\lambda_3 = \lambda_3$.
\\
If $a_2,a_3\geq 1$ then $a_2(-\lambda_2)+ a_3 (-\lambda_3) > -\lambda_3$. 
Note also that the cases $a_2=a_3=0$ and $a_2 =0,a_3=1$ are impossible.
\\
If $a_2=1$, $a_3 = 0$ then $\lambda_3 = a_1\lambda_1 + \lambda_2 \geq \lambda_2$, which is impossible by hypothesis. 

So for all resonance relations of type 3 we can take $i=2,l=1$, in view of Samovol's criterion. 

\medskip
We conclude using Theorem \ref{TheoremLinearizationSamavol} that $\bV$ is $C^1$ linearizable at $\overbp$.
\end{proof}
\section{Global boundary dynamics}
\label{sec:boundary dynamics}
As before, in this section we assume $k=2$ and explain the adjustments for general $k$ at the end of the section.
The rescaled geodesic vector field $\bV$ is tangent to the boundary $r=0$ of $\{E=\tfrac12\}\subset 
\threeT^*M$, which is the union of two components corresponding to $\xi=\pm1$.  We will 
study the dynamics of the flow of $\bV$ on the part $\xi=1$.
As explained in the introduction this can be identified with the space $T^*\partial M$ with variables $\vp, \theta$ in local coordinates.

Let 
\begin{equation}
\label{eqn:Ebounddef} 
 \Ebound (\vp,\theta)= S(\vp) + \tfrac12 \Cpartial^{-1}\theta^2  = (G_0)_{|\xi=1}.
\end{equation}

 By \eqref{eq: diffeq E1 leading terms}
the trajectories of the vector field $\bV$ on $T^*\dM$ are the solutions of the system
\begin{equation}\label{eqn:bd system}
 \vpdot = \Ebound_\theta, \qquad \thetadot = -\Ebound_\vp - 3\theta
\end{equation}
Explicitly,
\begin{equation}
\label{eqn:bd system explicit} 
 \vpdot = \Cpartial^{-1}\theta, \qquad \thetadot = -\Sbound_\vp (\vp)- 3\theta -\tfrac12 (\Cpartial^{-1})_{\vp}\theta^2.
\end{equation}
The system \eqref{eqn:bd system} is a Hamiltonian system modified by the damping term $-3\theta$. It 
describes motion of a particle on $\partial M$, with Riemannian metric given by $\gdM$, in the potential $S$ and in 
the presence of damping (for example, friction). 
The following properties are standard for such systems. We provide a proof for completeness. In this 
discussion the dimension of $M$ is arbitrary.
\begin{proposition}\label{prop:energy decay}\mbox{}
Let $M$ be a cuspidal manifold of order $k=2$. The singular points of $\bV$ on $T^*\partial M$ are equal to the critical 
points of the boundary energy $\Ebound$ and given by the critical points of $S$:
$$\{(\vp,\theta):\,\theta=0,\ S_\vp(\vp) = 0\}. $$
Every maximal trajectory $\gamma$ of $\bV$ on $T^*\partial M$ is defined for all times.
If $S$ is constant or has only isolated critical points then we have in addition:
\begin{enumerate}
 \item[a)]
 As $t\to\infty$, the trajectory 
 $\gamma(t)$ converges to a singular point of $\bV$. 
 \item[b)]
 As $t\to-\infty$, either $\Ebound$ is bounded along $\gamma$ and $\gamma(t)$ converges to a singular point of $\bV$,
 or else
 $\Ebound(\gamma(t))\to\infty$.
\end{enumerate}
\end{proposition}
% \ownremark{I (DG) removed the following from the statement since it is not needed in this paper: 
% If a non-constant trajectory starts at a singular point $(\vp,0)$ at $t=-\infty$ and ends at the singular 
% point $(\vp',0)$ at $t=\infty$ then $\Sbound(\vp') < \Sbound(\vp)$. 
% }
\begin{proof}
 Let $t\mapsto\gamma(t)=(\vp(t),\theta(t))$, $t\in(T_-,T_+)$ be a maximal trajectory of $\bV$. 
From \eqref{eqn:bd system} we get
\begin{equation}\label{eqn:lyapunov fcn}
   \Ebounddot = - 3|\theta|^2
\end{equation}
where $|\theta|^2 := \Cpartial^{-1}\theta^2$.
 In particular,  $\Ebound$ is decreasing along $\gamma$. Since $\Ebound$ is bounded below, it must approach 
a limit as $t\to T_+$. In particular, $\gamma$ is contained in a compact subset of $T^*\partial M$, which 
implies $T_+=\infty$. Fix $T\in(T_-,\infty)$.  Equation \eqref{eqn:lyapunov fcn} implies 
$\int_T^\infty |\theta(t)|^2\,dt<\infty$, and then the boundedness of $\thetadot$ easily implies that 
$\theta(t)\to0$ as $t\to\infty$. Let $\alpha(t)=S_\vp(\vp(t))$. The equation for $\thetadot$ and $\theta\to0$ 
imply that $\int_t^{t+1} \alpha(s)\,ds\to 0 $ as $t\to\infty$, and boundedness of $\dot{\alpha}$ again implies that 
$\alpha(t)\to 0$, so $\gamma$ converges to the set of critical points of $\Ebound$. If  this set is discrete 
then a) follows immediately. If $S$ is constant then $(\Cpartial^{-1})^{\cdot} = \vpdot (\Cpartial^{-1})_\vp = O(\theta)$, hence 
$$ (\Cpartial^{-1}\theta)^{\cdot} + 3\Cpartial^{-1}\theta = O(|\theta|^2) + \Cpartial^{-1}\thetadot + 3\Cpartial^{-1}\theta
= O(|\theta|^2) - \tfrac12 \Cpartial^{-1} (\Cpartial^{-1})_\vp \theta^2 = O(|\theta|^2)$$
which is integrable over $[T,\infty)$. Now $\theta\to0$ implies that $\lim_{T'\to\infty}\int_T^{T'} 
(\Cpartial^{-1}\theta)^{\cdot}\,dt $ exists, hence so does
$\lim_{T'\to\infty}\int_T^{T'} 3\Cpartial^{-1}\theta \,dt $ and hence  $\lim_{t\to\infty}\vp(t)$ exists. Choosing $T$ 
sufficiently big one may assume that all these curves lie in a single coordinate patch. This proves a).

Next, \eqref{eqn:lyapunov fcn} implies $\Ebounddot = -6\Ebound + O(1)$ which implies 
$|\Ebound(\gamma(t))|\leq K e^{6|t|}$ for some constant $K$. Then the same estimate must hold for $|\theta(t)|^2$, 
and this implies $T_-=-\infty$. Now by monotonicity either $\Ebound(\gamma(t))\to\infty$ as $t\to-\infty$, 
or else it is bounded. In the latter case \eqref{eqn:lyapunov fcn} implies that 
$\int_{-\infty}^\infty |\theta(t)|^2\,dt<\infty$. Then the same arguments as in  the proof of a) yield b).
\end{proof}
\begin{remark}\label{rem:perturbed bd system}
The proof carries over almost literally to the following perturbed version of b): Let $\gamma$ be a
solution of the system 
$$\vpdot = \Ebound_\theta + f(t),\quad \thetadot = -\Ebound_\vp - 3\theta - g(t)$$
where $|f(t)| + |g(t)| \leq Le^{-|t|/K}$ for all $t<0$, for some constants $K,L$. Then, if $\Ebound$ 
is bounded along $\gamma$, then $\gamma(t)$ converges to a singular point of $\bV$ as $t\to-\infty$.  
\end{remark}
% \ownremark{It is unclear to me (DG) whether this holds without the assumption of isolated critical points. 
% I would guess it does, but I am not sure. Of course $\gamma$ approaches the critical set, but in principle 
% it may happen that it comes closer and closer to it, with speed approaching zero but not converging to a 
% limit point.  One would have to show that $\vp(t)$ converges as $t\to\infty$, and the $\vpdot$ equation 
% shows that this would follow from the existence of $\int_T^\infty \theta(t)\,dt $. Unfortunately, this 
% does not follow from square integrability of $\theta$, but maybe it can be proven some other way?}
% \ownremark{One can get some more things out of this, for example (for surfaces): A trajectory starting from the smallest 
% maximum of $\Sbound$ must end up at a neighboring minimum (since it can never cross the $\vp$-value of 
% another maximum). Trajectories starting at a maximum can not make a full turn, i.e. $\vp$ cannot change 
% by length$(\partial\Xtilde)$. ('No spiraling, just wiggling') 
% }

We now analyze the case of surfaces $M$ in greater detail. Then we may and will, for simplicity, assume:
\begin{equation}\label{eqn:assn arclength}
 \vp \text{ is an arc length parameter for } \partial M. 
\end{equation}
 This means $\Cpartial\equiv 1$ on $\partial M$, so $\bV= \theta\partial_\vp - (S_\vp + 3\theta) \partial_\theta$ 
when restricted to ${\partial M}$.

%%%%%%%%%%%%%%%%%%%%%%%%%%%%%%%%%%%%%%%%%%%%%%%%%%%%
%%%%%%%%%%%%%%%%%%%%%%%%%%%%%%%%%%%%%%%%%%%%%%%%%%%%
\medskip
For the following discussion we assume:
\begin{equation}\label{eqn:assumptions}
\begin{gathered}
 \text{$S$ has a maximum at $\vp=0$ with $S_{\vp\vp}(0)<0$}\\
 \text{and a minimum at $\vpmin>0$, and $S_\vp<0$ on $(0,\vpmin)$.} 
\end{gathered}
\end{equation}
Denote  by $\bp_-, \bpmin$ the points in $\partial M$ given by $\vp=0$, $\vp=\vpmin$, respectively. 
Recall that these correspond to singular points $\overbp_-, \overbpmin$ of $\bV$ in $T^*\partial M$, 
given by these values of $\vp$, and  $\theta=0$. Recall that at a local maximum of $\Sbound$, the 
linear part of $\bV$ restricted to $T^*\partial M$ has eigenvalues $\lambda_3 < 0 < \lambda_2$, 
with eigenspace for $\lambda_2$ 
spanned by $\nu_2 = \partial_\vp + \lambda_2\partial_\theta$. Thus, by the unstable manifold theorem there 
are (up to time shift) precisely two trajectories starting at $\overbp_-$ at $t=-\infty$, 
and their directions at $\overbp_-$ are $\pm\nu_2$. We discuss only the trajectory starting in 
direction $\nu_2$ and denote it by $\Gammabpm$, there is an analogous 
discussion for the other trajectory.

\begin{figure}
 \includegraphics[width=6cm]{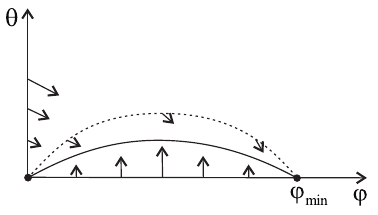}
 \caption{Boundary dynamics between a maximum of $S$ (at $\vp=0$) and a minimum in case $T=\infty$:
The dotted line is the graph of the barrier function $f$, the solid line is the trajectory $\Gammabpm$}
\label{fig:cusp boundary dynamics}
\end{figure}

The trajectory $t\mapsto\Gammabpm(t)=(\vp(t),\theta(t))$ starts out at $\overbp_-$ into 
the first quadrant $\vp>0$, $\theta>0$. It 
is defined for all times by Proposition \ref{prop:energy decay} a). As long as $\vp<\vpmin$, 
that is, for $t<T:= \sup\{t:\, \vp(s)<\vpmin\text{ for all }s\leq t\}$, it will stay 
in the first quadrant since $\bV$ points inside this quadrant along this part of its boundary, more 
precisely:
\begin{align*}
\bV (\vp,0) & =  -\Sbound_\vp (\vp)\partial_\theta \quad\text{ and } -\Sbound_\vp(\vp)>0 \text{ for } 0<\vp < \vpmin \\
\bV (0,\theta) & = \theta\partial_\vp -3\theta\partial_\theta \quad\text{and }\theta>0 \text{ for the first quadrant}.
\end{align*}
See Figure \ref{fig:cusp boundary dynamics}.
Since $\vpdot = \theta > 0$ in the first quadrant, $\vp$ is strictly increasing on $(-\infty,T)$. 
Now there are two possibilities: 
\begin{equation}
 \label{eqn: alternative}
\begin{aligned}
 \text{either } & T<\infty, \text{ then $\vp(T)=\vpmin$ and $\theta(T)>0$} \\
 \text{or }  &  T=\infty, \text{ then $\vp(t)<\vpmin\ \forall t$ and $\vp(t)\to\vpmin$, $\theta(t)\to 0$ as $t\to\infty$}
\end{aligned}
\end{equation}
(and hence $\Gammabpm(t)\to\overbpmin$), by Proposition \ref{prop:energy decay} a)
 since there are no critical points with $\vp\in(0,\vpmin)$.

We discuss the first case below in the proof of Theorem \ref{thm: S far from const}.
In the second case in \eqref{eqn: alternative} the  trajectory $\Gammabpm$ is the 
graph of a smooth function 
$h:(0,\vpmin)\to (0,\infty)$. 
We now give a sufficient condition for this to occur. It is obtained 
by a simple barrier argument.
\begin{lemma} \label{lem: barrier}
Suppose we are in the setup described above, starting from \eqref{eqn:assn arclength}. Assume that there is a $C^1$ 
function $f:[0,\vpmin] \to \R_+$ satisfying the conditions
\begin{gather} 
\label{eqn: cond a}
f(\vpmin)=0, \text{ and } f(\vp)>0 \text{ for } \vp\in (0,\vpmin)\\
 \label{eqn: cond c}
f(0)>0, \text{ or } f(0)=0 \text{ and } f'(0)>\lambda_2(\bp_-).\\
\label{eqn: cond b}
 ff' + 3f + \Sbound_\vp \geq 0\text{ on } (0,\vpmin)
\end{gather}
Then:
\begin{enumerate}
\item[a)]
 The closure of the image of $\Gammabpm$ is the graph of a $C^1$ function $h:[0,\vpmin]\to\R_+ $, 
and $h'(0)=\lambda_2(\bp_-)$. 
\item[b)]
If in addition $(ff' + 3f + \Sbound_\vp)' < 0 $ at $\vpmin$ then $\Gammabpm$ approaches $\bpmin$ 
tangent to the eigenvector $\nu_2(\bpmin)$, that is, $h'(\vpmin) = \lambda_2(\bpmin)$. 
\end{enumerate}

\end{lemma}
Concerning the  condition in b), observe that $f(\vpmin)=0$ implies differentiability of $ff'$ 
there, with derivative $(f'(\vpmin))^2$. Also, the function $g=ff'+3f+\Sbound_\vp$ satisfies $g(\vpmin)=0$ 
and $g(\vp)\geq 0$ for $\vp<\vpmin$, which implies that  $g'(\vpmin)\leq 0$. 
The given condition strengthens this to $g'(\vpmin)<0$. 
\begin{proof}
 Let
 $$ D(f) = \{(\vp,\theta): 0\leq \vp\leq \vpmin,\ 0\leq \theta \leq f(\vp)\} $$
 We claim that the vector field $\bV=\theta\partial_\vp -(\Sbound_\vp + 3\theta)\partial_\theta$ points 
inside $D(f)$, or is tangential to the boundary, at each boundary point of $D(f)$: This was checked 
before the proposition for the parts of the boundary where $\theta=0$ or $\vp=0$. The remaining part 
of the boundary is the graph of $f$, which is the zero set of $\ftilde(\vp,\theta)=f(\vp)-\theta$, 
with $\ftilde > 0 $ inside $D(f)$. Now $\bV \ftilde = \theta f' + \Sbound_\vp + 3\theta$, which at 
$\theta=f(\vp)$  equals $ff' + 3f + \Sbound_\vp$, so \eqref{eqn: cond b} gives $\bV\ftilde\geq 0$ which 
was to be shown. 
 
By \eqref{eqn: cond c} the trajectory $\Gammabpm$ starts out into the interior of $D(f)$, and 
by what we just proved it can never leave it. Since $\vpdot = \theta$ is positive in the interior of 
$D(f)$, the trajectory traces out the graph of a smooth function $\vp\mapsto h(\vp)$ on $(0,\vpmin)$. 
$h$ extends to be smooth at $0$ by the unstable manifold theorem, and it is $C^1$ at $\vpmin$ 
since the trajectory is tangent to an eigenspace at $\overbpmin$. This proves a).

Finally, $(ff' + 3f + \Sbound_\vp)'  = (f')^2 + 3f' + \Sbound_{\vp\vp}$ at $\vpmin$, and if this is negative 
then $f'$ (everything evaluated at $\vpmin$) must lie strictly between the zeroes $\lambda_2$ and $\lambda_3$ 
of the characteristic polynomial $\lambda^2 + 3\lambda + \Sbound_{\vp\vp}$ of the linear part of $\bV$ at 
$\vpmin$, so $h' \geq f' >\lambda_3$. Since $h'$ must be either $\lambda_2$ or $\lambda_3$ by the local 
behavior of trajectories, it follows that $h' = \lambda_2$ and hence b).
\end{proof}

\begin{proposition}\label{proposition-close-to-circle}
Let $M$ be a cuspidal surface of order $2$ and assume $S$ is a Morse function on $\partial M$ satisfying the bound 
$S_{\vp\vp}< 9/4$ everywhere, where $\vp$ is an arc length parameter on $\partial M$. Let
$\Gamma=\bigcup_\bp \Gamma_\bp$ where $\bp$ ranges over the critical points of $S$ and $\Gamma_\bp$ 
is the unstable manifold of the point $\overbp$ for the vector field $\bV$ on $T^*\partial M$. Then
$\Gamma$ is the graph of a $C^1$ section of $T^*\partial M$. That is,
$$ \Gamma = \{(\vp,\theta):\, \vp\in\partial M,\ \theta=h(\vp)\},\quad h\in C^1(\partial M).$$
Furthermore, $h'(\vp)=\lambda_2(\bp(\vp))$ at each critical point $\bp=\bp(\vp)$ of $S$.
\end{proposition}
In fact, the proof shows that $h$ is $C^\infty$ except at minima of $S$, and at the minima the regularity 
can, in general, be slightly improved, see the remark after Proposition \ref{prop:unstablemfdsmooth}.

% \ownremark{I (DG) removed some detail here to simplify the presentation: the statement can be localized 
% to parts of the boundary, and strict inequality $S{\vp\vp}<9/4$ is only needed at the minima of $S$.
% }
\begin{proof}
Since $S$ is a Morse function and $\dim\partial M=1$, its critical points are either maxima or minima, 
and these alternate along $\partial M$.
We first consider an interval between a maximum $\bp_-$ and the next minimum $\bpmin$. We may assume 
that the maximum is at $\vp=0$ and the minimum is at $\vpmin>0$. Then we are in the situation of Lemma 
\ref{lem: barrier}. As barrier function we take $f = -\tfrac23\Sbound_\vp$. This clearly satisfies 
condition \eqref{eqn: cond a}. Also, condition \eqref{eqn: cond b} holds since
$$ g:= ff'+3f+\Sbound_\vp = \tfrac49 \Sbound_{\vp}\Sbound_{\vp\vp} - 2\Sbound_{\vp} + \Sbound_{\vp}
= \Sbound_{\vp} (\tfrac49 \Sbound_{\vp\vp} - 1) > 0 $$
on $(0,\vpmin)$, and $g'(\vpmin)<0$ since $S_\vp=0$, $\tfrac49 \Sbound_{\vp\vp} - 1<0$, $S_{\vp\vp}>0$ at $\vpmin$.
Finally, since $f(0)=0$ condition \eqref{eqn: cond c} is equivalent to 
$-\tfrac23 a > \tfrac{-3+\sqrt{9-4a}}2$ where $a= \Sbound_{\vp\vp}(0)$, which is easily 
checked to be true for any $a<0$. So Lemma \ref{lem: barrier} a) and b) are applicable and yield a 
function $h_{\bp_-\bpmin}$ on the interval $[0,\vpmin]$ satisfying the claim of the Proposition there. 
We construct corresponding functions for each pair $(\bp_-,\bpmin)$ of consecutive maxima and minima of $S$, 
and by changing the signs also for pairs $(\bpmin,\bp_+)$ of consecutive minima and maxima; here 
$h_{\bpmin\bp_+}$ is negative in the interior of this interval. All these functions fit together to form a 
continuous function $h$ on all of $\partial M$ since they vanish at the endpoints of these intervals. $h$ 
is smooth at maxima $\bp$ of $S$ since its graph is the unstable manifold of $\overbp$ near $\overbp$, 
which is smooth. $h$ is $C^1$ at minima $\bpmin$ of $S$ since its derivatives from both sides agree, 
by Lemma \ref{lem: barrier}b).
\end{proof}
The value $\tfrac94$ is optimal for this type of barrier function $f$.
Note that the condition $0<\Sbound_{\vp\vp}<\frac94$ at a minimum $\bpmin$ is equivalent to the linear part of 
$\bV$ at $\bpmin$ having two distinct negative real eigenvalues.

See Section \ref{sec:cusp sing} for a sufficient condition when 
the condition of Proposition \ref{proposition-close-to-circle} is satisfied in the case of cuspidal singularities.

\subsection*{Adjustments for general order $k$}
The discussion generalizes to general $k$ by replacing $S$ with $\frac12 k(k-1)S$ and $3$ by $2k-1$, starting in \eqref{eqn:Ebounddef}, \eqref{eqn:bd system}. Compare \eqref{eqn:replace k}. Proposition \ref{prop:energy decay} generalizes literally (with obvious adjustments in the proof). In the statement of Lemma \ref{lem: barrier} replace $S$ and $3$ as described.

In Proposition \ref{proposition-close-to-circle} the bound $S_{\vp\vp}<9/4$ is replaced by $S_{\vp\vp}<a_k$, with $a_k$ defined in \eqref{eqn:def ak}. This results from replacing the barrier function $f=-\frac23 S_\vp$ by $f= - \frac {k(k-1)}{2k-1}S_\vp$.
\section{The exponential map in the case of a Morse function $S$}
\label{sec:exponential map}
In this section we prove Theorem \ref{thm: crit pts} for cuspidal manifolds of arbitrary dimension with $S$ a constant or a Morse function. Then we consider the case of a cuspidal surface, with $S$ a Morse function. We define the exponential map and prove 
Theorems \ref{thm:surfaces}, \ref{thm: S close to const} and \ref{thm: S far from const}.

\subsection*{Proof of Theorem \ref{thm: crit pts}}
Assume $k=2$ first. The geodesics are, up to reparametrization, the projections to $M$ of the integral curves of the vector 
field $\bV$ on the energy hypersurface $\{E=\tfrac12\}\subset\threeT^*M$. Thus, let $\gamma:(T_-',T_+')\to\threeT^* M$ be a 
maximal integral curve of $\bV$, and assume that $\gamma_{|(T_-',T')}$ is contained in $\{r<r_0\}$ for 
some $T'\in(T_-',T_+')$. The number $r_0$ will be chosen sufficiently small below. Then necessarily $T_-'=-\infty$ 
since $\bV$ is tangent to the boundary. We will 
show that $\lim_{t\to-\infty}\gamma(t)$ exists and is a singular point of $\bV$. 

Write $\gamma(t)=(r(t),\vp(t),\xi(t),\theta(t))$. Here we have chosen a trivialization of $M$ near the boundary, 
which allows us to identify $\threeT^*M$ with $T^*[0,r_0)\times T^*\partial M$ there, and then $(r,\xi)$ denotes a 
point in the first factor and $(\vp,\theta)$ a point in the second factor. The proof falls into two parts: First, we use the 
$\rdot$ and $\xidot$ equations in \eqref{eq: diffeq E1 leading terms} to show that $\xi(t)$ must be 
close to one on $(-\infty,T')$, hence $r(t)$ exponentially decreasing as $t$ decreases, and that $\theta(t)$ 
is bounded. Then the decay of $r$ together with the boundedness of $\theta$ shows that the $\vpdot$, $\thetadot$ 
dynamics is close to the boundary dynamics, which allows us to use Proposition \ref{prop:energy decay}b).

By \eqref{eq: diffeq E1 leading terms}  there is $K>0$ such that $\xidot \geq 2(1-\xi^2) - Kr$ 
on the energy hypersurface. 
\smallskip

\noindent Claim 1: $\xi>\frac12$ for $t<T'$.
\begin{proof}
For $r_0$ sufficiently small we have  $\xidot>1$  for $|\xi|\leq\frac12$, so if $\xi(t_0)\leq \frac12$ 
for some $t_0<T'$ then $\xi(t)<-\frac12$ for all $t<t_0-1$, and then $\rdot=\xi r + O(r^3)$ implies 
that $r$ leaves the interval $[0,r_0]$ as $t\to-\infty$, if $r_0$ is chosen sufficiently small (depending on the constant in $O(r^3)$). This would 
 contradict the assumption that $ r(t)<r_0$ for all $t<T'$.
 \end{proof}

\smallskip

\noindent Claim 2: $\xi^2 > 1 - Kr$ for $t<T'$.
\begin{proof}
 Let $f(r,\xi)=\xi^2 + Kr$. The part of $\gamma_{|(-\infty,T')}$ where $f>1$ is forward invariant under 
the flow of $\bV$, since at any point with $f(r,\xi)=1$, we have
$2(1-\xi^2)-Kr=Kr$ and therefore 
$ \bV f =2\xi\xidot + K\rdot\geq 
2\xi(2(1-\xi^2) -  Kr) + K(r \xi + O(r^3)) 
= 3\xi Kr + O(r^3) > 0 $, 
% $\frac12 \bV f =\xi\xidot + K\rdot\geq \xi Kr + K(r \xi + O(r^3)) \geq \
% Kr + O(r^3) > 0 $, 
by Claim 1, for $r_0$ sufficiently small. 
Hence if $f\leq 1$ at some point $\gamma(t_0)$ 
then this remains true for all $t\leq t_0$. Now $f(r,\xi)\leq 1$ implies $2(1-\xi^2)-Kr \geq 1-\xi^2>0$, 
so we would get $\xidot > 1-\xi^2>0$ for all $t\leq t_0$. Then $\xi$ would be monotone increasing in 
$t\leq t_0$, so $\lim_{t\to-\infty}\xi(t)$ would exist and lie in $[\frac12,1)$. However, 
this is impossible since it would imply $\xi(t_0)-\xi(-\infty)>\int_{-\infty}^{t_0}(1-\xi(t)^2)dt = \infty$ and 
hence unboundedness of $\xi$.
\end{proof}

Now we conclude from \eqref{eqn: energy shell} and Claim 2 that
$r\theta^2$ is bounded (everything for $t<T'$). Using this in the $\dot\xi$ equation in \eqref{eq: diffeq E1 leading terms} we see that in fact $\dot\xi \geq 2(1-\xi^2)-Kr^2$, and then essentially the same argument as in Claim 2 shows that we even have
$\xi^2>1-Kr^2$. This in turn, again in conjunction with \eqref{eqn: energy shell}, shows that $\theta$ is bounded.
From Claim 1 and $\rdot=r\xi+O(r^3)$ we obtain $r(t)\leq Ke^{-|t|/3}$. Also, \eqref{eqn: energy shell} implies that $|\xi(t) - 1| \leq K'e^{-2|t|/3}$. 
Therefore Remark \ref{rem:perturbed bd system}, applied to the $\vpdot$ and $\thetadot$ equations in 
\eqref{eq: diffeq E1 leading terms}, implies that $\gamma(t)$ approaches a singular point of $\bV$ 
as $t\to-\infty$. This proves parts (a) and (b) of Theorem \ref{thm: crit pts}.

The claim that each singular point $\overbp$ of $\bV$ is the starting point of a trajectory of 
$\bV$ follows from the invariant manifold theorem and Lemma \ref{lemma-eigenvalues} since the 
unstable tangent space at $\overbp$ contains the eigenvector for the eigenvalue $\lambda=1$, 
which points inside $M$.

Finally, returning to the unrescaled geodesic vector field we have $r'=\xi + O(r^2)$ along a geodesic, 
and since $\xi$ is near one this implies that the geodesic reaches $r=0$ in finite backward time, so $T_-$
is finite. Also, it implies that there is $\tau_0>0$ only depending on $r_0$ so that $r(t)<r_0$ for $t-T_-\leq \tau_0$, and this proves part (c).

This completes the proof of Theorem \ref{thm: crit pts} in the case $k=2$. 
For general $k$ the proof is analogous, using \eqref{eqn: energy shell k}, \eqref{eq: diffeq E1 leading terms k} instead of \eqref{eqn: energy shell}, \eqref{eq: diffeq E1 leading terms}.
\subsection*{Definition of the exponential map}
We now define the exponential map in the case where $M$ is a cuspidal surface (of any order $k\geq2$) with connected boundary $\dM$ 
and $S$ is a Morse function on $\dM$. For this we need to parametrize appropriately 
the family of geodesics starting at the points of $\dM$. Recall that geodesics start at the 
critical points of the function $S$, and are projections of integral curves  of the rescaled geodesic vector field $\bV$ on the energy level $\{E=\frac{1}{2}\}\subset\threeT^*M$ (resp. ${}^{(2k-1)}T^*M$). 

The boundary $\dM$ is diffeomorphic to a circle. Fixing an orientation of $\dM$, we label the critical points 
of $S$ in cyclic order $\bp_1,\bp_1',\bp_2,\bp_2',\dots,\bp_l,\bp_l'$, where $l\in\N$ and 
where $S$ has maxima at the $\bp_i$ and minima at the $\bp_i'$. 

Recall the considerations in the proof of Proposition \ref{prop-localmax}, 
in the case of surfaces: Each critical point $\bp$ of $S$ is the projection of a hyperbolic singular point $\overbp$ of $\bV$. The geodesics hitting the boundary $\dM$ in $\bp$ are, as sets, the projections of the integral curves of $\bV$ whose union form
the interior of the unstable manifold $M^u_\bp$ of $\bV$ at $\overbp$. This unstable manifold is two dimensional for each maximum $\bp=\bp_i$
and projects diffeomorphically to $M$ in a neighborhood of $\overbp_i$. 
It is one-dimensional for each minimum $\bp=\bp_i'$.

Recall from Theorem \ref{thm: crit pts} that there is $\tau_0>0$ so that all geodesics starting at $\partial M$ exist at least for time $\tau_0$.
In the following definition we consider the circle  $\bS^1$ as the interval $[0,2\pi]$ with the endpoints identified.
\begin{definition} \label{def:exp map}
Under the assumptions and with the notation introduced above we define the exponential map as follows:
\begin{enumerate}
 \item
 For each maximum $\bp_i$ of $S$ we parametrize the set of geodesics starting at $\bp_i$ by the interval
 $I_i = \left(\frac{2(i-1)\pi}l,\frac{2i\pi}l\right)\subset \bS^1$, preserving orientation.
 \item
 We let the unique geodesic starting at the minimum $\bp_i'$ of $S$ correspond to the point $\frac{2i\pi}l\in S^1$.
 \item
 This defines a parametrization $q\mapsto \gamma_q$, $q\in \bS^1$, of the geodesics starting at $\dM$. Then we define
 $$ \exp_{\dM}:\bS^1\times(0,\tau_0) \to U\setminus \dM, \quad (q,t) \mapsto \gamma_q(\tau) $$
 where $\gamma_q$ is parametrized by arclength, starting at $\dM$ at time $\tau=0$.
\end{enumerate}

\end{definition}
Note that we did not specify the precise parametrization in (1). The reason for this is that there is no natural parametrization. In fact, even the choice of intervals $I_i$ is arbitrary (apart from their ordering on the circle). More precisely, if $\kappa:\bS^1\to \bS^1$ is any orientation preserving homeomorphism, then $\exp'_\dM(q,\tau) = \exp_\dM(\kappa(q),\tau)$ will serve just as well as exponential map. The set of exponential maps $\exp_\dM(q,\tau)=\gamma_q(\tau)$ obtained in this way are characterized by the 

\begin{quote}
 {\bf Local order preserving property:}
  If $a,b,c\in \bS^1$ are pairwise different and lie in this order on $\bS^1$ then there is $\tau_1>0$ so that the geodesics
$\gamma_a,\gamma_b,\gamma_c$ do not intersect and lie in this order for times in $(0,\tau_1)$. 
\end{quote}
Here, $\tau_1$ may depend on 
$a,b,c$, and for certain cuspidal metrics cannot be chosen uniformly for all $a,b,c$. This leads to 
discontinuity of the exponential map, see Remark \ref{rem:exp discontinuous}.

In Remark \ref{rem:exp blowup} we sketch a more intrinsic definition of the exponential map, which makes use of a generalized notion of inhomogeneous blow-up.
\smallskip

To clarify part (1) of the definition, we now give an explicit parametrization using linearizations.

Let  $\bp = \bp_i$ be a maximum of $S$, and let $\Utilde\subset M_\bp^u$ be a neighborhood of $\overbp$ so that the projection $\pi:\Utilde\to M$ is a diffeomorphism to its image. The rescaled geodesic vector field $\bV$ is tangent to $\Utilde$. Let $\bV_{|\Utilde}$ be its restriction to $\Utilde$ and $\bV' = \pi_*(\bV_{|\Utilde})$ its projection to $\pi(\Utilde)\subset M$.
Since the linearization of $\bV$ at $\overbp$ has eigenvalues $1$ and $\lambda_2(\bp)>0>\lambda_3(\bp)$, the linearization of $\bV_{|\Utilde}$ at $\overbp$ has eigenvalues $1$ and $\lambda_2:=\lambda_2(\bp)$. Then the same is true for the linearization of $\bV'$ at $\bp$.
Since both eigenvalus are positive, $\bV'$ is $C^1$-linearizable near $\bp$  by Hartman's Theorem \cite{Har} (see also \cite[p 127]{Per}). 

Therefore, we may choose $C^1$ coordinates $\rho\geq0,\psi$ near $\bp$, where $\rho=0$ on $\partial M$, so that in these coordinates
$\bp=(0,0)$ and 
$\bV' = \rho \dpl_\rho + \lambda_2 \psi\dpl_\psi$.
We may also assume that the  coordinates are chosen so that the orientation of $\dM$ 
corresponds to the positive $\psi$-direction.
The integral curves (as sets) of $\bV'$ not contained in the boundary $\rho=0$ 
are of the form $\{\psi = K \rho^{\lambda_2},\rho >0\}$ with $K\in \R$.
Therefore, $K$ parametrizes the family of integral curves of $\bV'$, hence of geodesic 
leaving $\bp$, and does so in an order preserving way. Now choosing a diffeomorphism $\R\to I_i$, 
for example $K \mapsto \frac{\pi}l\left(\frac K{1+K^2} + 2i-1\right)$, we get a parametrization 
as required in (1) of Definition \ref{def:exp map}.
\begin{remark}\label{rem:exp blowup}
 The domain of the parametrization $q\mapsto\gamma_q$ may be described somewhat more naturally with 
the help of the following construction: For a vector field $\bV$ on a smooth manifold $\Sigma$ with an 
unstable critical point $\bp$ there is a natural notion of blow-up of $\bp$ in $\Sigma$ with respect to 
$\bV$. This is a smooth manifold with boundary, denoted $[\Sigma,\bp]_\bV$, together with a blow-down 
map $\beta_\bp:[\Sigma,\bp]_\bV\to\Sigma$. The boundary (front face) parametrizes the integral curves of $\bV$ starting 
at $\bp$. It is diffeomorphic to a sphere and may be thought of as small sphere around $\bp$ transversal 
to $\bV$. The map $\beta_\bp$ may not be smooth, but $\bV$ lifts to a smooth vector field on 
$[\Sigma,\bp]_\bV$, given by $s\partial_s$ near the boundary for a suitable boundary defining 
function $s$. The construction generalizes to the case where $\Sigma$ is a manifold with boundary, 
$\bp\in\partial\Sigma$, and $\bV$ is tangent to $\partial\Sigma$, and then yields a manifold 
with corners. See \cite[Section 2]{HMV} for details.

We apply this to the surface $M^u_{\bp_i}$ for each $i$ and obtain $[M^u_{\bp_i},\overbp_i]_\bV$. 
The front face $F_i:=\beta_\bp^{-1}(\overbp_i)$ of the blow-up is diffeomorphic to a closed interval, and is oriented by the orientation of $\dM$. We now glue the right endpoint of each $F_{i-1}$ to the left endpoint of $F_i$ (with $F_0=F_l$) and obtain a manifold $F$ homeomorphic to $\bS^1$ which is the natural domain of parametrization of geodesics leaving $\dM$. Then the exponential map is defined on $F\times (0,\tau_0)$.

The exponential map in the sense of Definition \ref{def:exp map} is obtained by parametrizing the interior of $F_i$ by the interval $I_i$ for each $i$.

%From a different perspective, we may blow up each $\bp_i$ in $M$ with respect to the projection of the vector field $\bV$ near $\bp_i$. Then $T$ is naturally identified with the space obtained by collapsing all intervals on $\dM$ on adjacent front faces. 

\end{remark}

\subsection*{Proof of Theorem \ref{thm:surfaces}}
We constructed the exponential map above. We need to show that it is surjective. For this we prove 
that the union of the unstable manifolds $M^u_\bp$ over all critical points projects {\em onto} a neighborhood 
of $\partial M$ under the canonical projection $\pi:\threeT^*M\to M$ (resp. $\pi:{}^{(2k-1)}T^*M\to M$ for general $k$).

First, observe that each unstable manifold $M_\bp^u$ intersects the boundary $\{r=0\}$ transversally 
near $\overbp$ since this is true for the linear part of $\bV$, and since both $M_\bp^u$ and $\{r=0\}$ 
are invariant under the flow, the intersection is transversal everywhere.
The intersection $M_\bp^u\cap \{r=0\}$ is the unstable manifold $\Gamma^u_\bp$ of $\overbp$ of $\bV$ 
restricted to the boundary.

The projections of the unstable manifolds $M^u_{\bp_i'}$ are (closures of) single geodesics starting 
at $\bp_i'$, which divide a neighborhood of $\dM$ into connected components $U_i$, each containing 
one maximum $\bp_i$. It suffices to show that $\pi(M^u_{\bp_i})$ contains $U_i$ for each $i$. 

For this, first consider the intersections with the boundary. 
Fix $i$ and write $\bp_-=\bp_i$ and $\bpmin=\bp_i'$. By the discussion before Lemma \ref{lem: barrier}, 
the part  $\Gamma_{\bp_-}$ of $\Gamma^u_{\bp_-}$ leaving $\bp_-$ in the positive $\vp$ direction projects 
onto the (open) interval from $\bp_-$ to $\bpmin$. A similar statement holds for the interval from 
$\bp_{i-1}'$ to $\bp_-$. Therefore,  $M^u_{\bp_i}\cap\{r=0\}$ projects onto $U_i\cap\{r=0\}$. Since 
$M^u_{\bp_i}$ intersects the boundary transversally, it projects onto some neighborhood of the open 
interval from $\bp_{i-1}'$ to $\bp_i'$. We need to show that this neighborhood cannot shrink to zero 
width when one approaches $\bp_{i-1}'$ or $\bp_i'$. We will prove that $U_i$ intersected with a 
neighborhood of $\bp_i'$ is contained in $\pi(M^u_{\bp_i})$; the argument at $\bp_{i-1}'$ is then analogous.

Recall the dichotomy \eqref{eqn: alternative}. If $T<\infty$ here then $\pi(\Gamma_{\bp_-})$ 
contains $\bp_i'$, so we are done. On the other hand, if $T=\infty$ then $\Gamma_{\bp_-}(t)\to\overbpmin$ 
as $t\to\infty$, and the eigenvalues $\lambda_2(\bpmin)$, $\lambda_3(\bpmin)$ are real. In this case, 
the behavior of $M^u_{\bp_-}$ near $\overbpmin$ may be understood by linearizing $\bV$ near $\overbpmin$. 
We use Lemma \ref{lem:inv mfd smooth} below, applied as explained after its statement. For later purposes 
the lemma is stated more strongly than needed here. Here we only need the consequence that, in a neighborhood 
of $\overbpmin$, the curve $M^u_\bpmin$ is contained in the closure of $M^u_{\bp_-}$, so the projection 
of $M^u_\bpmin$ is contained in the closure of $U_i$. This completes the proof of Theorem \ref{thm:surfaces}.

%%%%%%%%%%%%%%%%%%%%%%%%%%%%%

\begin{lemma}
 \label{lem:inv mfd smooth}
 Let $\lambda, \mu<0 $ and let $V$ be the linear vector field 
$x_1 \partial_{x_1} + \lambda x_2 \partial_{x_2} +\mu x_3 \partial_{x_3}$ 
on $\R^3$. Let $m\geq 1$. Let $\Sigma$ be an invariant $C^m$ surface with boundary 
$\Gamma=\Sigma\cap\{x_1=0\}$, and assume this intersection is transversal and $\Gamma$ is a single 
trajectory of $V$. Denote by $X_1$ the non-negative $x_1$-axis.
 
Then a neighborhood $\Sigma'$ of zero in $\Sigma \cup X_1$ is a surface with corner. More precisely, let  
$\tau = -\frac{1}{\lambda}$ and $\rho = \frac{\mu}{\lambda}$, so $\tau,\rho>0$. Assume 
$\Gamma\subset\{x_2>0\}$. Then $\Sigma'$ is a graph
\begin{equation}\label{eqn:Sigma as graph} 
\begin{gathered}
  \Sigma' = \{(x_1,x_2,x_3):\ x_3 = A(x_1x_2^\tau) x_2^\rho 
,\ (x_1,x_2)\in U'\} 
\end{gathered}
\end{equation}
where  $A$ is a $C^m$ function on $[0,\varepsilon)$ for some $\varepsilon>0$, and $U'$ is a 
neighborhood of the origin in the quarter plane $\R_+\times \R_+$.
\end{lemma}
The lemma is applied as follows, in the context of the discussion before Lemma \ref{lem: barrier}: 
Let $\bpmin$ be a minimum of $S$, and suppose the second alternative in \eqref{eqn: alternative} holds. 
Then $\Gamma_{\bp_-}(t)$ approaches $\overbpmin$ as $t\to\infty$. By Proposition \ref{linearization-loc-min}
we may linearize $\bV$ near $\overbpmin$ by a $C^m$-diffeomorphism where $m\in\{1,2\}$. Thus, we introduce local coordinates 
$x_1,x_2,x_3$ so that $\overbpmin$ is the origin and $\bV$ is given by 
$x_1 \partial_{x_1} + \lambda x_2 \partial_{x_2} +\mu x_3 \partial_{x_3}$ on a neighborhood $U$ where 
$\{\lambda,\mu\} = \{\lambda_2(\bpmin),\lambda_3(\bpmin)\}$ and the  boundary is given by $\{x_1=0\}$.
We take $\Gamma=\Gamma_{\bp_-}$ and $\Sigma=M^u_{\bp_-}$. Clearly $M^u_{\bpmin}=X_1$. 
Since $\Gamma$ is a trajectory of $\bV$ approaching the origin, it must have a tangent vector there, 
and we may assume that the coordinates are chosen so that this vector is $\partial_{x_2}$. Then 
$\Gamma\subset\{x_2>0\}$ near the origin, so the assumptions of the lemma are satisfied. 
\begin{proof}[Proof of Lemma \ref{lem:inv mfd smooth}]
Choose a point $p$ of $\Gamma$. Since the intersection of $\Sigma$ with $\{x_1=0\}$ is transversal we 
may choose a $C^m$-curve $\omega(s)=(s,\omega_2(s),\omega_3(s))$, $s\in I=[0,\eps)$ contained in $\Sigma$, 
with $\omega(0)=p$, so $\omega_2(0)>0$. Also, since $x_1x_2^\tau$ is constant along any integral curve 
of $\bV$, only trajectories passing through this curve will contribute to $\Sigma'$ if this is chosen 
sufficiently small.

For $s\in I$ let $\R_+\to \R^3$, $t\mapsto\Gamma_s(t)$ be the forward integral curve of $\bV$ starting 
at $\omega(s)$, so $\Gamma_s(0)=\omega(s)$ and $\Gamma_0=\Gamma$, up to time shift.

For each $s$, the quantities $x_1 x_2^{\tau}$ and $x_3 x_2^{-\rho}$ are constant along $\Gamma_s$, that is, 
for points $(x_1,x_2,x_3)$ in the image of $\Gamma_s$. Evaluating at $t=0$ shows that, along $\Gamma_s$, 
\begin{equation}\label{eqn:integrals}
 x_1 x_2^\tau = s \omega_2(s)^\tau,\quad x_3 x_2^{-\rho} = \omega_3(s) \omega_2(s)^{-\rho}. 
\end{equation}
Now $\omega_2(0)>0$, hence the equation $\sigma = s \omega_2(s)^\tau$ can 
be solved $C^m$-smoothly for $s=s(\sigma)$, for $\sigma$ in a half neighborhood $J$ of zero, and the function
$A:J\to \R$, $\sigma \mapsto \omega_3(s(\sigma))\omega_2(s(\sigma))^{-\rho}$ is $C^m$. Then we have along 
$\Gamma_s$, with this $\sigma$, $$ x_3 = A(\sigma)x_2^\rho = A(x_1 x_2^\tau) x_2^\rho.$$
The claim follows.
\end{proof}

%%%%%%%%%%%%%%%%%%%%%%%

\subsection*{Proof of Theorem \ref{thm: S close to const}}
Assume $k=2$ first. The heart of the proof is that the various unstable manifolds fit together nicely.
So we first prove the following proposition.

\begin{figure}
 \includegraphics[width=6cm]{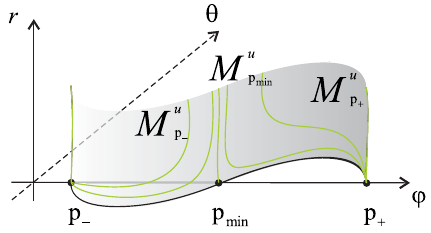}%Dynamics-near-boundary-Nr_2.eps}
%  %% here it didn't work to include the .eps since the shaded regions made the printers stall
 \caption{The dynamics near the boundary.}
 \label{fig:cusp dynamics}
\end{figure}

\begin{proposition}\label{prop:unstablemfdsmooth}
Assume the setting and the conditions of Theorem \ref{thm: S close to const}.
 Let $M^u$ be the union of the unstable manifolds $M_\bp^u$, where $\bp$ ranges over the critical points 
of $S$. Then $M^u$ is a $C^1$ manifold, and there is $r_0>0$ so that the 
intersection of $M^u$ with $\{r<r_0\}$ projects diffeomorphically to $M\cap\{r<r_0\}$ under the 
projection $\pi:\threeT^*M\to M$.
\end{proposition}
See Figure \ref{fig:cusp dynamics}.

The proof actually gives more regularity than $C^1$: The manifold $M^u$ is $C^\infty$ except at 
points on the curves $M_\bp^u$ where $\bp$ is a minimum of $S$, 
and here the regularity is  $C^{\alpha(\bp)}$ where
\begin{equation}
 \label{eqn:regularity}
 \alpha(\bp) = 
\begin{cases}
 1 & \text{ if } \lambda_2\in \Q \\
 \min\{2,\frac{\lambda_3}{\lambda_2}\} & \text{ if } \lambda_2 \in \R\setminus\Q
\end{cases}
\end{equation}
with $\lambda_2=\lambda_2(\bp)$, $\lambda_3=\lambda_3(\bp)$.
%; recall that $-\frac32 < \lambda_2 < 0$ at a non-degenerate minimum of $S$
% \ownremark{note that the case $m=\rho=2$ does not occur since $\rho=2$ implies $\lambda_2=-1$ and hence $m=1$}
%
%
\begin{proof}[Proof of Proposition \ref{prop:unstablemfdsmooth}]

The main task is to prove that the different $M_\bp^u$ fit together to form a $C^1$ manifold, 
the projection statement will then be seen to be an easy consequence.

Since integral curves do not intersect, the manifolds $M_\bp^u$ do not intersect for different $\bp$. 
Since $M^u$ is invariant under the flow, it is enough to prove the regularity near $r=0$. 

First consider the boundaries $\Gamma_\bp^u=\partial M_\bp^u= M_\bp^u\cap\{r=0\}$ and their union 
$\Gamma^u = \bigcup_\bp \Gamma_\bp^u$. These are the unstable manifolds of $\bV$ restricted to the
boundary $T^*\partial M$. By assumption $S_{\vp\vp}< \frac94$ everywhere, so we may apply Proposition 
\ref{proposition-close-to-circle}  and conclude that $\Gamma$ is the graph of a $C^1$ function 
$\Gamma=\{(\vp,h(\vp)):\, \vp\in\partial M\}$. Thus, both the regularity and the projection statement 
of Proposition \ref{prop:unstablemfdsmooth} are true for the boundary of $M^u$.

We now prove that $M^u$ is a $C^1$ manifold  in a neighborhood of $r=0$, with boundary $\Gamma^u$. 
Since $M_\bp^u$ is a smooth surface for each maximum $\bp$ of $S$, we only need to prove the regularity 
in a small neighborhood $U$ of $\overbpmin$ where $\bpmin$ is a minimum of $S$. Thus, let $\bpmin$ be a 
minimum of $S$ and let $\bp_-$, $\bp_+$ be the maxima of $S$ closest to $\bpmin$, so that the corresponding 
parameter values satisfy $\vp_-<\vpmin<\vp_+$. Then the set $M^u \cap U$ is the union of $M_{\bp_-}^u 
\cap U$, $M_{\bpmin}^u \cap U$ and $M_{\bp_+}^u \cap U$. We now apply Lemma \ref{lem:inv mfd smooth} as 
explained after its statement, first to $\Sigma_-=M_{\bp_-}^u\cap U$ and then to $\Sigma_+=M_{\bp_+}^u\cap U$. 
By the last statement of Proposition \ref{proposition-close-to-circle}, $\Gamma_{\bp_\pm}$ have tangents 
$\pm\nu_2(\bpmin)$ at $\overbpmin$, so we take $\lambda=\lambda_2(\bpmin)$, $\mu = \lambda_3(\bpmin)$ in 
Lemma \ref{lem:inv mfd smooth}, and this implies $\rho>1$ there. The lemma shows that $\Sigma_\pm$ are graphs 
as in \eqref{eqn:Sigma as graph}, with possibly different functions $A_\pm$. 

The more precise regularity statement \eqref{eqn:regularity} follows from the regularity in Lemma 
\ref{lem:inv mfd smooth} together with the regularity of the linearization, Proposition \ref{linearization-loc-min}.

Finally we prove that $M^u$ projects diffeomorphically to $M$ near $r=0$. For any boundary point 
$\overbq\in\Gamma^u=\dM^u$, the tangent space $T_\overbq M^u$ is transversal to $\{r=0\}$, since this is 
already true for each $M_\bp^u$ separately. In addition, $T_\overbq M^u \supset T_\overbq \Gamma^u$,
so $T_\overbq M^u\cap T_\overbq \{r=0\}$ contains a vector with nonzero $\partial_\vp$-component. These 
two facts combine to show that, for each $\overbq\in\dM^u$, the tangent space  $T_\overbq M^u$ projects 
isomorphically to $T_\bq M$ under $d\pi$, where $\pi:\threeT^*M\to M$ is the projection, and 
this implies the claim.
\end{proof}

We now prove Theorem \ref{thm: S close to const}.
By Theorem \ref{thm: crit pts} every geodesic starting at $\partial M$ must do so at a critical point 
$\bp$ of $S$, so the corresponding integral curve of $\bV$  starts at a singular point $\overbp$ of $\bV$.
This integral curve is then contained in the unstable manifold $M_\bp^u$ of $\overbp$. Conversely, 
$M_\bp^u\cap \{r>0\}$ is  the union of such integral curves. Therefore, Proposition \ref{prop:unstablemfdsmooth} 
implies that $\expdM$ is bijective. 

It remains to prove the continuity of $\expdM$ and of its inverse. We continue to work on $M^u$. Some care 
needs to be taken since the vector field $\bW$ blows up at the boundary. But the main point is that only 
its $\partial_\vp$ component blows up, while its $\partial_r$ component is $\xi+O(r^2)=1+O(r^2)$. Recall that 
the domain of $\expdM$ is the set of $(q,\tau)$ where $\tau\in(0,\tau_0)$ and $q$ is a point on a 
transversal $T_i$ near $\overbp_i$ for some $i$, where the transversals are glued at their endpoints. 
The continuity of $\expdM$ at points $(q_0,\tau)$ where $q_0$ is not a boundary point of a 
transversal $T_i$ is clear, so we assume that $q_0$ is a boundary point of $T_i$, say the 'right' 
boundary point, which labels the geodesic starting at the next minimum $\bp_i'$. The idea is this: 
As sets, the geodesics $\gamma_q$, with $q$ in the interior of $T_i$ approaching $q_0$, converge to 
the union of the boundary trajectory $\Gamma_{\bp_i}$ and of $M^u_{\bp_i'}$. Since $r'=\frac{dr}{d\tau}$ 
is approximately one, the part of $\gamma_q$ near $\Gamma_{\bp_i}$ is traversed in a very short time. 
Then the part near $M^u_{\bp_i'}$ must, including its time parametrization, be close to $\gamma_{q_0}$.

More precisely, fix $(q_0,\tau)$ and a neighborhood $U$ of $\gamma_{q_0}(\tau)$. Denote by 
$r$, $r_{0}$ the $r$-components  of integral curves $\gamma$, $\gamma_{q_0}$ of $\bW$.
 Then there are $\delta>0$, a neighborhood $U'$ of $M^u_{\bp_i'}$  and a neighborhood $V'$ of 
$\tau$ so that for  all $0<\tau_1<\tau$ and all integral curves $\gamma$ of $\bW$ lying in $U'$ 
for the time interval $(\tau_1,\tau)$ and satisfying $|r(\tau_1) - r_0(\tau_1)|<\delta$ we have that 
$\gamma(\tau')\in U$ for all $\tau'\in V'$. This is because $r'=F(r,\vp,\theta)$ where $F$ is $C^1$ 
and the values of $\vp,\theta$ for $\gamma$ and $\gamma_0$ will be close together at all times in 
$(\tau_1,\tau)$ if $U'$ is chosen sufficiently small.

Next, for the $\delta>0$ and neighborhood $U'$ of $M^u_{\bp_i'}$ obtained above there is a neighborhood 
$V$ of $q_0$ in $T_i$ so that for all $q\in V$ the curve $\gamma_q$ first runs inside $\{r<\delta/2\}$ 
and then inside $U'$. Now $\frac{dr}{d\tau} = 1+O(r^2)$ implies that the travel time $\tau_1(q)$ for the 
first part is at most on the order of $\delta/2$. Then the $r$-components $r_q$ and $r_{0}$ of $\gamma_q$, 
$\gamma_{q_0}$ satisfy $|r_q(\tau_1(q)) - r_0(\tau_1(q))| < \delta$, so the first part of the argument can
be applied with $\gamma=\gamma_q$ and $\tau_1=\tau_1(q)$.
Summarizing, given any neighborhood $U$ of $\gamma_0(\tau)$ we have found neighborhoods $V$ 
of $q_0$ and $V'$ of $\tau$ so that $\gamma_q(\tau')\in U$ for all $q\in V$, $\tau'\in V'$. This 
proves the continuity of $\expdM$, from the left with respect to $q$. Continuity from the right follows 
by the same argument applied to the left endpoint of $T_{i+1}$. The continuity of the inverse is proved 
in a similar way, but easier. For example, the second component, $\tau$, of $\expdM^{-1}$ is simply the 
distance to the boundary, and its continuity follows from the triangle inequality.  

This concludes the proof of Theorem \ref{thm: S close to const} in the case $k=2$. The argument for general $k$ is exactly the same, except that one replaces $\threeT^*M$ by ${}^{(2k-1)}T^*M$, and in the proof of Proposition \ref{prop:unstablemfdsmooth} one uses
the condition  $S_{\vp\vp}<a_k$ instead of $S_{\vp\vp}<9/4$ when referring to Proposition \ref{proposition-close-to-circle}.

\begin{remark}\label{rem:exp discontinuous}
The proof  also shows why the exponential map may be discontinuous in general: Suppose $l\geq 2$ in Definition \ref{def:exp map}, 
i.e. the function $S$ has at least two maxima and at least two minima. Suppose $\Gamma_{\bp_1}$ 
does not approach $\bp_1'$ but rather continues in the upper half plane $\theta>0$ and then 
approaches $\bp_2'$. Examples of functions $S$ yielding this boundary dynamics can easily be 
constructed. Let $q_1,q_2$ be the labels of the geodesics leaving $\bp_1'$, $\bp_2'$, respectively. 
Then, for any fixed $\tau>0$, the point $\gamma_q(\tau)$ will approach $\gamma_{q_2}(\tau)$ rather 
than $\gamma_{q_1}(\tau)$ as $q\to q_1$ from the left, so $\expdM$ is discontinuous at $(q_1,\tau)$ for any $\tau>0$.
Also, it is easy to see that this discontinuity cannot be removed by a simple reordering (i.e. 
by a different gluing prescription), which in any case would be unnatural since it would break 
up the order preserving property of $\expdM$.
\end{remark}

\subsection*{Proof of Theorem \ref{thm: S far from const}}
Let $\bpmin$ be a point where $S$ has a local minimum and where $S_{\vp\vp}>a_k$. Then the eigenvalues 
$\lambda_2(\bpmin)$, $\lambda_3(\bpmin)$ are non-real with real part $-\frac{2k-1}2$, so nearby trajectories of 
$\bV$ on $T^*\partial M$ spiral towards $\overbpmin$.

Consider the maxima $\bp_-$, $\bp_+$ of $S$ closest to $\bpmin$, where $\vp_-<\vpmin<\vp_+$ for 
the corresponding parameters. W.l.o.g.\ we may assume $S(\vp_-)\leq S(\vp_+)$. Then we are in the 
situation discussed before Lemma \ref{lem: barrier} (where we put $\vp_-=0$ for simplicity), more
precisely in the first case of \eqref{eqn: alternative}. We claim that the trajectory $\Gammabpm$ 
considered there approaches $\overbpmin$ as $t\to\infty$. To show this, recall from the proof of Proposition 
\ref{prop:energy decay} that the boundary energy function $\frac12 k(k-1)S(\vp) + \theta^2/2$ and hence $S(\vp)+\frac1{k(k-1)}\theta^2$ is strictly decreasing along 
$\Gammabpm$. At $t=-\infty$ this function has the value $S(\vp_-)$, and on the line $\vp=\vp_+$ 
its values are at least $S(\vp_+)\geq S(\vp_-)$. Hence, $\Gammabpm$ can not cross or approach this 
line, and neither the line $\vp=\vp_-$. So the only possible limit point of $\Gammabpm$ for $t\to\infty$ 
is $\overbpmin$, which was to be shown.

Now note that the spiraling of $\Gammabpm$ around $\overbpmin$ as $t\to\infty$ implies that 
the projection of $\Gammabpm$ to the $\vp$-axis oscillates around $\vpmin$ infinitely often.

Now consider the unstable manifold $M_{\bp_-}^u$. The image of $\Gammabpm$ is an open subset of 
$\partial M_{\bp_-}^u$. By continuity of $\bV$, for any $T>0$ and $\eps>0$ there are interior 
trajectories of $\bV$ contained in $M_{\bp_-}^u$ that follow the boundary trajectory $\Gammabpm$ 
within an error of $\eps$ on the time interval $[0,T]$. In particular, for any $\eps>0$ there is a 
trajectory $\gamma$ starting at $\overbp_-$ whose projection to $M$ intersects the projection of 
the trajectory $M^u_\bpmin$ at a point with $r<\eps$. This proves the theorem.
%
% \ownremark{in general, in the first case of \eqref{eqn: alternative} the trajectory may approach 
% another singular point of $\bV$.
% A sufficient condition for the first case in \eqref {eqn: alternative} to occur is that the linear part of 
% $\bV$ at $\overbp_+$ has non-real eigenvalues, since then any trajectory approaching $\overbp_+$ must spiral 
% around this point. 
% }
%
%
%
%
%
%
%
%
%
%
%
%%%%%%%%%%%%%%%%%% &&&&
%
%
%
%
%
%
%
%
%
\section{Application to cuspidal singularities}\label{sec:cusp sing}
The main motivation for defining cuspidal metrics is that they arise from cuspidal singularities.
The notion of cuspidal singularity itself does not involve a metric. 
In this section we give two definitions of cuspidal singularity, prove their equivalence and show that the restriction of a smooth ambient metric to a cuspidal singularity yields
a cuspidal metric upon resolution of the singularity, so the main theorems are applicable in this setting. Also, we interpret the quantities $S$ and $\bg_{\partial M}$ for cuspidal singularities and give a sufficient condition for the exponential map to be a homeomorphism.

\subsection{Definition of cuspidal singularities}
We now discuss the notion of cuspidal singularity for subsets $X$ of a smooth manifold $Z$. Since this is a local notion, one may always think of $Z=\R^n$; however, it is useful to have an invariant geometric understanding.

Recall that the tangent cone $C(X,p)$ of a subset $X\subset \R^n$ at $p\in X$ is the set of limits of secant half-lines 
through $p$ from points in $X$:
\begin{center}
\vspace{4pt}
$C (X,p) := \{ r\nu:\, r>0, \nu\in \Sr^{n-1}, \exists (p_j)_j \in X\setminus\{p\}: p_j\to p$ and
$\lim_j \frac{p_j-p}{|p_j -p|} =\nu\}$.
\vspace{4pt}
\end{center}
It is well-known and easy to check that for subsets $X$ of a manifold $Z$ the tangent cone at $p\in X$ is invariantly defined as a subset $C(X,p)\subset T_pZ$. 
\begin{definition}
\label{def:cusp sing}
Let $k\geq2$, $k\in\N$. 
Let $X$ be a subset of a manifold $Z$ and $p\in X$.  We say that $X$ has a 
{\em cuspidal singularity of order $k$} at $p$ if 
\begin{enumerate}
 \item[(a)]
the tangent cone $C(X,p)$ is a 
half-line and 
 \item[(b)]  there are coordinates $(x_1,\dots,x_{n-1},z)$ on a neighborhood $U\subset Z$ of $p$ in which $p=0$ and $C(X,p)=\R_{>0}\partial_z$, so that $X\cap U$ is given locally as in
\eqref{eqn:beta}, \eqref{eqn:beta'}.
\end{enumerate} 
\end{definition}
Here $\partial_z=(0,\dots,0,1)$.
If these conditions are satisfied then we say that
the singularity of $X$ at $p$ can be resolved by a {\em blow-up of order $k$}. Note that condition (b) implies that $(X\cap U)\setminus \{p\}$ is a submanifold of $Z$.
Below we need another characterization of cuspidal singularities in terms of a standard (or first order) blow-up followed by a blow-up of order $k-1$, which also has the virtue of leading to a manifestly invariant characterization.  To formulate this we recall some basic terminology.

\subsubsection{Review of blow-up}
The (oriented, first order) {\em blow-up} of a manifold $Z$ in a point $q\in Z$ is a geometric, coordinate free counterpart to introducing polar coordinates around $q$. By definition it is a manifold with boundary, denoted $[Z,q]$, together with 
a smooth map $\beta_q: [Z,q]\to Z$ (the {\em blow-down map}) which maps the boundary $\partial[Z,q]$ to $q$ and 
is a diffeomorphism from $[Z,q]\setminus \partial [Z,q]$ to $Z\setminus \{q\}$, and which is, locally near 
$\partial[Z,q]$ resp. $q$, given by the following model:  $Z=\R^n$, $q=0$, and then  $[Z,q]=\R_{+}\times \bS^{n-1}$ 
and $\beta_q (r,\omega) = r\omega.$ Thus $\partial[Z,q]=\{0\}\times \bS^{n-1}$ in the model. See \cite{mel:aps} or \cite{gri:BBC} 
for a more in-depth discussion, in particular of the coordinate invariance of this notion and of its generalization
to manifolds with corners. Of this we only need the case where $Z$ is itself a manifold with boundary and 
$q\in\partial Z$. Then $[Z,q]$ is a manifold with corners; the local model is $Z=\R^{n-1}\times \R_+$, $q=0$, 
$[Z,q]=\R_{+}\times \bS^{n-1}_+$, where $\bS^{n-1}_+=\bS^{n-1}\cap(\R^{n-1}\times\R_+)$ is the upper half sphere
and $\beta_q$ is as before. In this case, $\beta_q$ maps $\{0\}\times \bS^{n-1}_+$ to $0$ and is a diffeomorphism
between the complements of these sets. In either case $\beta_q^{-1}(q)$ is called the {\em front face} of the blow-up.

Starting from a coordinate system for $Z$ near $q$, {\em projective coordinates} are defined for $[Z,q]$ as follows.
The coordinates for $Z$ identify a neighborhood of $q$ with $\R^n$, with $q=0$.
Denote the coordinates by $(x,z)$ where $x=(x_1,\dots,x_{n-1})$. 
Then the {projective coordinates} 
$(y,z)$, $y=(y_1,\dots,y_{n-1})\in\R^{n-1}$, $z\in\R_+$ are defined to be the unique coordinates on the 'upper half' $U= 
\{(r,\omega)\in\R_+\times \bS^{n-1}:\, \omega_n>0\}$ of $[\R^n,0]$ in terms of which
 the blow-down map is $\beta_q(y,z)=(zy,z)$. (So informally $y_j=\frac{x_j}z$ for $j=1,\dots,n-1$. The relation to the 
$r,\omega$-variables is $y_i=\frac r{\omega_n}\omega_i$, $z=r\omega_n$, but this is not needed.) There 
are also projective coordinate systems covering the remainder of $[\R^n,0]$, but we do not need them here.

The blow-up to order $k$ can be given a similar invariant description, denoted $\beta_q:[Z,q]_k\to Z$,
and \eqref{eqn:beta} is the blow-down map in order $k$ projective coordinates.

If $X\subset Z$ then the {\em strict transform} of $X$ under the blow-up of $q$ is the closure of the 
pre-image of $X\setminus \{q\}$, so $\beta_q^*X := \overline{\beta_q^{-1}(X\setminus\{q\})}$. A 
{\em p-submanifold} (p for 'product type') of a manifold with boundary $Z$ is a submanifold $X\subset Z$ such that 
$\partial X\subset\partial Z$ and $X$ hits $\partial Z$ transversally. This extends to manifolds with 
corners $Z$; we only need the straightforward case where $X$ intersects the boundary only in the interior 
of a boundary hypersurface.

\subsubsection{Alternative characterization of cuspidal singularities}
Let $Z$ be a manifold and $X\subset Z$ a subset. We say that $X$ has a {\em conical singularity} at $q\in X$ if it is resolved by blowing up $q$. That is,
there is an open neighborhood $U\subset Z$ of $q$ so that the strict transform of $X\cap U$  is a p-submanifold of $[U,q]$. If $Z$ is a manifold with
boundary and $q\in\partial Z$ then we also require non-tangency of $X$ to $\partial Z$, that is, that the boundary of the strict transform of $X\cap U$ be contained in the interior 
of the front face of $[U,q]$.

\begin{lemma} \label{lem:quadr blow-up}
Let $X$ be a subset of a manifold $Z$ and $p\in X$. 
%Let $k\geq 2$, $k\in\N$.
Let $Z'=[Z,p]$ be the blow-up of $Z$ in $p$ with blow-down map $\beta_p:Z'\to Z$. 
Then $X$ has a cuspidal singularity of order $k$ at $p$ if and only if 
 its strict transform 
$X'=\beta_p^*(X)\subset Z'$ intersects $\partial Z'$ in a single point $q$ and 
\begin{list}{}{}
\item if $k=2$: $Z'$ has a conical singularity at $q$
\item if $k>2$: $Z'$ has a cuspidal singularity of order $k-1$  
 at $q$.
\end{list}
\end{lemma}
Thus, the singularity of $X$ can be resolved by first blowing up $p$ and then blowing up $q$ to order $k-1$. The resolution is defined as
the strict transform 
\begin{equation}
 \label{eqn:resolution}
\Xtilde=(\beta_p\circ\beta_{q})^*X .
\end{equation}
By iteration, this implies that an order $k$ blow-up can be replaced by a sequence of $k$ standard blow-ups. We will not use this consequence.
\begin{proof}
Since this is a local statement we may assume $Z=\R^n$, $p=0$. 
Denote $X'=\beta_p^*(X) \subset Z'=[\R^n,p]$.
Points of the front face $\partial Z'$ correspond to 
directions at $p$, so $C(X,p)= \{r\nu:\, \nu\in X'\cap \partial Z',\ r>0\}$. Therefore, condition (a) of Definition \ref{def:cusp sing} is equivalent to $X'$ intersecting $\partial Z'$ in a single point $q$.  Assuming this is satisfied we may choose coordinates $(x,z)$, $x=(x_1,\dots,x_{n-1})$ on $\R^n$ so that  $C(X,p)=\R_{>0} \partial_z$. In the corresponding  projective coordinates $(y,z)$, for which $\beta_p(y,z)=(zy,z)$, this corresponds to $q=0$. 

The coordinates $(y,z)$ on $Z'$ define projective coordinates $(u,z)$ on $[Z',q]_{k-1}$, for which the blow-down map $[Z',q]\to Z'$ is given by $\beta_{q}(u,z)=(z^{k-1}u,z)$, hence $\beta_p(\beta_q(u,z)) = (z^ku,z)$. 
Let $\tilde X\subset[Z',q]_{k-1}$ be the strict transform of $X'$ under the order $k-1$ blow-up of the manifold with boundary $Z'$ in $q$.  The coordinates $(u,z)$ are defined in a neighborhood of the interior of the front face of $[Z',q]_{k-1}$. Therefore, $X'$ having a conical singularity at $q$ implies that a neighborhood of the boundary of $\tilde X$  is contained in the domain of definition of these coordinates. 

This shows that condition (b) of Definition \ref{def:cusp sing} is simply the condition that $X'$ has a conical (resp.\ order $k-1$ cuspidal) singularity at $q$, written in the coordinates $(u,z)$. The map $\beta$ in \eqref{eqn:beta} is the map $\beta_p\circ\beta_q$.
\end{proof}

\subsection{The relation between cuspidal singularities and cuspidal metrics}
The motivation for introducing the notion of cuspidal metric is the following proposition.

\begin{proposition} \label{prop:cuspidal sing}
Let $X$ be a subset of a manifold $Z$ having a cuspidal singularity of order $k$ at $p\in X$. Assume $X\setminus \{p\}$ is a submanifold of $Z$. 
Let $\Xtilde$ be the manifold with boundary obtained as resolution of $X$ as in \eqref{eqn:resolution}. 

Let $g$ be a smooth Riemannian 
metric on $Z$ and $g_X$ its restriction to $X\setminus\{p\}$. Then the pullback $g_\Xtilde$ of $g_X$ 
to $\Xtilde\setminus\partial\Xtilde$ extends to a cuspidal metric of order $k$ on $\Xtilde$.
\end{proposition}
In other words, after resolving the singularity one can choose
coordinates $(r,\varphi)$ near any boundary point of the resolution $\Xtilde$ so that the metric takes the form \eqref{eqn:def-cusp-metric}.

\begin{proof}
We use the characterization of a cuspidal singularity of order $k$ given in Lemma \ref{lem:quadr blow-up}. We proceed as follows: We use  geodesic polar coordinates on the first blow-up $Z'=[Z,p]$, then projective coordinates on the second blow-up $\Ztilde=[Z',q]_{k-1}$, then modify these to make the mixed terms of the metric vanish to order $2k$. In the end we restrict to $\Xtilde\subset\Ztilde$.

Since everything is local near $p$, we may assume $Z=\R^n$, $p=0$. 
Denote the blow-down maps $\beta_p:Z'\to\R^n$, $\beta_{q}:\Ztilde\to Z'$,
$\beta=\beta_p\circ\beta_q$ as in \eqref{eqn:resolution}. We need to find a suitable $r$-coordinate on $\Ztilde$, defined near the interior of its front face which is then given by $r=0$, 
so $\beta^*g$ has the form \eqref{eqn:def-cusp-metric}.
 
Introduce geodesic polar coordinates for $g$ near $0$. This means that we choose a diffeomorphism of 
$Z'=[\R^n,0]$ with $\R_+\times \bS^{n-1}$ so that, with $R$ the coordinate on $\R_+$, the metric has the form
\begin{equation}
\label{eqn:normal coords} 
 \beta_p^*g = dR^2 + R^2 h(R)
\end{equation}
where $R\mapsto h(R)$ is a smooth family of Riemannian metrics on $\bS^{n-1}$ (here $dR^2 :=(dR)^2$ as usual). In fact, $h(0)$ is the 
standard metric on $\bS^{n-1}$, but this is inessential for the construction.

Next let $U=(U_1,\dots,U_{n-1})$ be any coordinates, centered at $q$, on  $\bS^{n-1}$. Then
$h(R) = \sum_{ij} h_{ij}(R,U) dU_i \, dU_j$ with $h_{ij}$ smooth. On $\Ztilde=[Z',q]_{k-1}$ use 
projective coordinates $R,u_i=\frac{U_i}{R^{k-1}}$; they are defined in a neighborhood of the interior 
of the front face.%, so in particular on $\Xtilde$ near $\partial \Xtilde$.

Then $(\beta_q^* h_{ij})(R,u)=h_{ij}(R,R^{k-1}u)=h_{ij}(0,0)+O(R)$ and $\beta_q^*dU_i = d(R^{k-1}u_i) = u_i (k-1)R^{k-2} \,dR + R^{k-1}\,du_i$ and so, with 
\begin{equation}
\label{eqn:def S h'}
 S(u)=\sum_{ij} h_{ij}(0,0) u_i u_j,\quad h'=\sum_{ij} h_{ij}(0,0) du_i du_j,
\end{equation}
a short calculation yields
\begin{equation}
\label{eq:h metric}
\beta_q^*h(R) =R^{2k-4} \left[(k-1)^2 S \, dR^2 + (k-1)R\,dR dS + R^2 h'
+ \htilde \right]
\end{equation}
where $\htilde = O(R) (dR)^2 + \sum_i O(R^2)\,dR du_i + \sum_{ij} O(R^3)\, du_i du_j $.  Then from \eqref{eqn:normal coords} we have 
\begin{equation}
 \label{eqn:xxx}
\begin{aligned}
\beta^*g&=\beta_q^*\beta_p^*g
= dR^2 + R^2 \beta_q^*h(R) \\
&
= (1+(k-1)^2 R^{2k-2}S)\, dR^2 + (k-1)R^{2k-1}\,dR dS + R^{2k}\,h' + R^{2k-2} \tilde h\,.
\end{aligned}
\end{equation}
Since the $\htilde$ part only contributes $O(R^{2k-1})dR^2 + \sum_i O(R^{2k})dRdu_i + \sum_{ij}O(R^{2k+1})du_idu_j$ this has the form \eqref{eqn:def-cusp-metric}, except for the term $(k-1)R^{2k-1} \,dRdS$ (and with different $R^{2k-2} dR^2$ coefficient). 
We would like to get rid of this term by setting $r=R + wR^{2k-1}$. Equation  \eqref{eqn:metric change} 
shows that for this we need to take $w=\frac{k-1}2 S$, and then $-k(k-1)S + 2(2k-1)w = (k-1)^2S$ is precisely the coefficient of $R^{2k-2}dR^2$ in \eqref{eqn:xxx}.
Thus we obtain \eqref{eqn:def-cusp-metric}, where the coefficient of $r^{2k-2} dr^2$ is now $-k(k-1)S$ and $c_{ij}(r,u) = h_{ij}(R,R^{k-1}u)$.

Now consider the submanifold $\Xtilde\subset \Ztilde$. Since it is transversal to $r=0$, it can be locally
parametrized as $u=u(r,\vp)$, $\vp=(\vp_1,\dots,\vp_{m-1})$, with 
$\vp\mapsto u(0,\vp)$ having injective differential. Here $m=\dim \tilde X$.
%the matrix $\left(\frac{\partial u_i}{\partial\varphi_j}\right)_{i=1,\dots,n-1; j=1,\dots,m-1}$ having full rank $m-1$ at $r=0$. 
Restricting $\beta^*g$ to $\Xtilde$ amounts to writing $(\beta^*S)(r,\vp) = S(r,u(r,\vp))$ and similarly for the other coefficient functions, and $du_i =
\sum_{j=1}^{m-1} \frac{\partial u_i}{\partial\varphi_j} d\varphi_j$ and
therefore yields a cuspidal metric again.
Invariantly, $S$ is simply restricted to $\Xtilde$ and the metric $\sum_{ij}c_{ij} d\vp_i d\vp_j$ on $\partial \Xtilde$ is the restriction of the metric $\sum_{ij}c_{ij} du_i du_j$ on the interior of the front face of $\Ztilde$.
\end{proof}
\begin{remark}\label{rem:S, gdm for Xtilde}
The proof shows that for a cuspidal manifold $\Xtilde$ arising as resolution of a space $X$ with cuspidal singularity $p$ the quantities 
$S$ and $\gdXtilde$ (see Lemma \ref{lem:cusp metric}) are given as follows. First assume that the ambient space is $Z=\R^n$ with $g$ the standard 
Euclidean metric, that $p=0$ and $C(X,p)=\R_{>0}\partial_z$ in coordinates $(x,z)$ as in Definition \ref{def:cusp sing}(b).
Let $X_z = \{x\in\R^{n-1}: (x,z)\in X\}$ be the cross section of $X$ at height $z$. 
Then 
$\partial\Xtilde = \lim_{z\to 0+} z^{-k}X_z\subset \R^{n-1}$ in the Hausdorff sense, and 
\begin{equation}
\label{eqn:S gdXtilde embedded} 
 S(\bu) = |\bu|^2,\quad \gdXtilde = \text{ restriction of }\geucl\text{ to } \partial\Xtilde 
 \end{equation}
where $\geucl$, $|\ |$ are the Euclidean metric and norm on $\R^{n-1}$.

More generally and invariantly, the ambient metric $g$    induces the structure of Euclidean vector space on the front face of the 
order $k$ (or iterated standard) blow-up of the ambient space, and then $S$ and $\gdXtilde$ are given by \eqref{eqn:S gdXtilde embedded} for that Euclidean metric.

The dotted lines in Figure \ref{fig:cusp} indicate the tangent cone of $X$ at $p$ and the line $\bu=0$ in its resolution.
\end{remark}
\begin{remark}
 The proof shows that Proposition \ref{prop:cuspidal sing} holds for the more general class of metrics obtained from any 
conical metric by blowing up a boundary point to order $k-1$.
\end{remark}

\subsection{Cuspidal singularities with convex base}
% \ownremark{Next, $S_\vp = 2u\cdot u_\vp$, so if $S_\vp=0$ at $\bp$ then $N_\bp=\frac {u(\bp)}{|u(\bp)|}$ is a unit 
% normal to the curve $\partildX$ at $\bp$. If $n=3$, where $\partildX$ is a plane curve, 
% $\kappa_\bp = -N_\bp\cdot H_\bp$ is the curvature of $\partildX$ at $\bp$ with respect to the normal 
% pointing towards the origin, so $\kappa>0$ if the curve bends toward the origin, and we obtain
% $$ 1 - ac = |u(\bp)| \kappa_\bp $$
% For example, if $|u|\equiv R$ then $\kappa=R^{-1}$ everywhere, so $1-ac=1$ which is clear from 
% $a=\frac12 S_{\vp\vp}=0$. }

Recall from Remark \ref{rem:S, gdm for Xtilde} that $\partial\Xtilde$ is naturally the subset of a 
Euclidean vector space. If $X$ is given as in \eqref{eqn:beta} then this is simply $\R^{n-1}\times\{0\}$ 
with the standard Euclidean metric. In general we still may identify it with $\R^{n-1}$ with the standard metric.
\begin{theorem}\label{thm:convex bound}
Let $X$ be a surface with cuspidal singularity satisfying the following assumptions.
\begin{enumerate}
 \item
 $\partildX\subset \R^{n-1}$ is contained 
in the boundary of a strictly convex set which contains the origin. 
\item
For any point $\bp\in\partial\Xtilde$ where $\partial\Xtilde$ is tangent to the 
sphere through $\bp$ centered at the origin, it is only simply tangent to that sphere.
\end{enumerate}
Then the exponential map based at the cuspidal singularity $q$ of $X$ and associated with any ambient metric restricted to $X$  is a homeomorphism near $q$.
\end{theorem}
In the case of a surface in $\R^3$ the simple tangency condition is equivalent to the condition that the 
osculating circle of $\partial\Xtilde$, wherever it is defined, be never centered at the origin. Strict convexity is meant in the sense of nonzero curvature.
\begin{proof}
Choose an arc length parametrization $\vp\mapsto \bu(\vp)$ of $\partial\Xtilde$. Then $S(\vp)=|\bu(\vp)|^2$ by Remark \ref{rem:S, gdm for Xtilde}, 
so the simple tangency condition is equivalent to $S$ being a Morse function.
Next,
\begin{equation}
 \label{eqn:meaning Spp}
\tfrac12 S_{\vp\vp}= |\bu_\vp|^2 + \bu\cdot \bu_{\vp\vp} = 1 +  \bu \cdot K(\bu)
\end{equation}
 for  $\bu=\bu(\vp)$ where $K(\bu)$ is 
the curvature vector of $\partial\Xtilde$  at $\bu$.
Suppose $\partildX$ is contained in the boundary of the strictly convex set $\calK$, and let $P_\bu$ be a supporting 
hyperplane for $\calK$ through $\bu$.  Then $K(\bu)$ points into that closed half-space determined by $P_\bu$ 
which contains $\calK$. Since $\calK$ contains the origin, it follows that $\bu \cdot K(\bu) < 0$, so 
$\Sbound_{\vp\vp}(\vp) < 2 < a_k$ for all $\vp$, where $k$ is the order of the cuspidal singularity. The claim now follows from Theorem \ref{thm: S close to const}. 
\end{proof}
The proof shows that the conclusion also holds with a certain amount of non-convexity since $a_k>2$. However, the larger $k$ the sharper the convexity assumption becomes, because $a_k\to2$ as $k\to\infty$.
\section{Examples}
\label{sec:examples}

We consider  surfaces $X\subset\R^3$ given as in \eqref{eqn:beta} with
$\Xtilde=\partildX\times\R_+$ a cylinder, for different boundary curves 
$\partildX\subset\R^2$. We use the Euclidean metric on $\R^3$. Write coordinates on $\R^2$ as 
$u=(v,w)$. Recall that $S(v,w) = v^2 + w^2$, see Remark \ref{rem:S, gdm for Xtilde}.
\begin{enumerate}
 \item
 If $\partildX$ is a circle centered at the origin then $S$ is constant, the geodesics
 of $X$ starting at the origin foliate the surface. This is obvious by rotational symmetry. 
 \item
If $\partildX$ is an ellipse centered at the origin, 
\begin{center}
\vspace{4pt}
$\displaystyle{\frac{v^2}{c^2} + \frac{w^2}{b^2}= 1}$, with $c>b> 0$.
\vspace{4pt}
\end{center}
The function $S$ has two maxima at $(\pm c,0)$ and two minima at $(0,\pm b)$. The boundary $\partildX$ 
is simply tangent to the circles centered at the origin passing through these points. Since $\partildX$ 
bounds a strictly convex set containing the origin, the exponential map based at $0$ is a local homeomorphism
by Theorem \ref{thm:convex bound}. There is one geodesic starting from each of the points $(0,\pm b)$, 
all the others start at $(\pm c,0)$.
\item
We now consider a circle not having the origin in its interior,
$$ (v-c)^2 + w^2 = 1,\quad c>1.$$
The function $S$ has a  minimum at $\bpmin=(c-1,0)$ and a maximum at $\bpmax=(c+1,0)$. These are 
non-degenerate, so we have simple tangency again. In the arc length parametrization 
$v(\vp)=c + \cos\vp$, $w(\vp)=\sin\vp$, we find $S(\vp) = c^2+1+ 2c\cos\vp $, so $S_{\vp\vp}(\vp)= -2c\cos\vp$.
 At the minimum $\vp=\pi$ of $S$ we see that $S_{\vp\vp}=2c$, which is also the maximal value of $S_{\vp\vp}$.
 We obtain the following picture:

There is a single geodesic $\gammamin$ starting at $\bpmin$, all others start from $\bpmax$. 
Their behavior depends on $c$:
\\
If $1<c<9/8$ % and $2c \neq 2$ can be removed since follows from c>1
 then Theorem \ref{thm: S close to const} is applicable. 
\\
If $c>9/8$ then we are in the situation of Theorem \ref{thm: S far from const}. Geodesics starting 
at $\bpmax$ almost tangentially to $\partildX$ will approach $\bpmin$, then oscillate around 
$\gammamin$ many times before they escape far away from the singularity.
\end{enumerate}

\smallskip
Our first and simplest example above might suggest that Theorem \ref{thm:S=const} is not so interesting. However, the point of the theorem is that all that matters for the conclusion is the constancy of $S$ {\em on the boundary}. Hence, any surface $\Xtilde$ with the same boundary will have a foliation by geodesics near the singularity. As another example, consider a surface $\tildX$ embedded in $\R^n$ with $n>3$, with $\partildX$ any curve in a sphere $R\cdot \bS^{n-2}$. Then there is no rotational symmetry, not even for $\partildX$, but Theorem \ref{thm:S=const} still gives a foliation by geodesics.
%
%
%
%
%
%
%
%
%
%
%
%
%
%
%
%
%
%
%
%
%
%
%
%
%
%
%
%
%
%
%
%
%
%
%

%\bibliographystyle{plain}
%\bibliography{dglib,mypapers}
%\end{document}

\end{document}